\documentclass[a4paper, 11pt, reqno]{amsart}
\usepackage{amsthm,amssymb,amsmath,enumerate,graphicx,psfrag,enumitem,mathrsfs}

\usepackage[margin=2.8cm]{geometry}
\usepackage{hyperref}
\hypersetup{colorlinks=true,
   citecolor=blue,
   filecolor=blue,
   linkcolor=blue,
   urlcolor=blue
}
\usepackage{dsfont}
\usepackage[textsize=footnotesize]{todonotes}

\usepackage{algorithm}
\usepackage{tikz}
\usetikzlibrary{arrows,patterns,decorations.markings}
\usetikzlibrary{backgrounds}

\numberwithin{equation}{section}

\newtheorem{definition}{Definition}[section]
\newtheorem{theorem}[definition]{Theorem}
\newtheorem{prop}[definition]{Proposition}
\newtheorem{lemma}[definition]{Lemma}
\newtheorem{cor}[definition]{Corollary}
\newtheorem{fact}[definition]{Fact}
\newtheorem{problem}[definition]{Problem}
\newtheorem{conj}[definition]{Conjecture}

\theoremstyle{definition}
\newtheorem{defin}[definition]{Definition}

\newtheoremstyle{claimstyle}{5pt}{5pt}{\em}{5pt}{\em}{:}{5pt}{}
\theoremstyle{claimstyle}
\newtheorem{claim}{Claim}

\newtheoremstyle{stepstyle}{10pt}{5pt}{\em}{0pt}{\em}{:}{5pt}{}
\theoremstyle{stepstyle}

\newcommand{\comment}[1]{}
\newcommand{\N}{\mathbb N}

\newcommand{\Z}{\mathbb Z}

\newcommand{\cG}{\mathcal{G}}
\newcommand{\cP}{\mathcal{P}}

\newcommand{\cF}{\mathcal{F}}

\newcommand{\cC}{\mathcal{C}}

\newcommand{\cH}{\mathcal{H}}
\newcommand{\cW}{\mathcal{W}}
\newcommand{\cM}{\mathcal{M}}
\newcommand{\cI}{\mathcal{I}}

\newcommand{\eul}{{\rm e}}

\newcommand{\cA}{\mathcal{A}}
\newcommand{\cU}{\mathcal{U}}
\newcommand{\cX}{\mathcal{X}}

\newcommand{\fa}{\mathfrak{a}}
\newcommand{\fc}{\mathfrak{c}}

\newcommand{\cupdot}{\mathbin{\mathaccent\cdot\cup}}

\newcommand{\defect}{{\rm def}}

\newcommand{\prob}[1]{\mathrm{\mathbb{P}}\left[#1\right]}
\newcommand{\expn}[1]{\mathrm{\mathbb{E}}\left[#1\right]}

\def\eps{{\varepsilon}}

\def\sm{\setminus}
\newcommand{\Set}[1]{\{#1\}}
\newcommand{\set}[2]{\{#1\,:\;#2\}}

\def\In{\subseteq}
\newcommand{\vv}[1]{{\vec{#1}}}

\renewcommand{\epsilon}{\varepsilon}

\def\noproof{{\unskip\nobreak\hfill\penalty50\hskip2em\hbox{}\nobreak\hfill%
       $\square$\parfillskip=0pt\finalhyphendemerits=0\par}\goodbreak}
\def\endproof{\noproof\bigskip}

\def\noclaimproof{{\unskip\nobreak\hfill\penalty50\hskip2em\hbox{}\nobreak\hfill%
       $-$\parfillskip=0pt\finalhyphendemerits=0\par}\goodbreak}

\def\endclaimproof{\noclaimproof\medskip}

\newdimen\margin
\def\textno#1&#2\par{
   \margin=\hsize
   \advance\margin by -4\parindent
          \setbox1=\hbox{\sl#1}
   \ifdim\wd1 < \margin
      $$\box1\eqno#2$$
   \else
      \bigbreak
      \hbox to \hsize{\indent$\vcenter{\advance\hsize by -3\parindent
      \it\noindent#1}\hfil#2$}
      \bigbreak
   \fi}

\def\proof{\removelastskip\penalty55\medskip\noindent\setcounter{claim}{0}\setcounter{step}{0}{\bf Proof. }} 

\def\lateproof#1{\removelastskip\penalty55\medskip\noindent\setcounter{claim}{0}\setcounter{step}{0}{\bf Proof of #1. }} 

\def\claimproof{\removelastskip\penalty55\medskip\noindent{\em Proof of claim: }}

\newcommand{\defn}{\emph}
\newcommand{\sub}{\subseteq}

\newcommand{\COMMENT}[1]{}

\newcounter{stepc}

\title{Resolution of the Oberwolfach problem}

\author[S.~Glock]{Stefan Glock}
\address[S.~Glock]{Institute for Theoretical Studies, ETH Z\"urich, Switzerland}
\email{stefan.glock@eth-its.ethz.ch}

\author[F.~Joos]{Felix Joos}
\address[F.~Joos]{Institut f\"ur Informatik, Universit\"at Heidelberg, Germany }
\email{joos@informatik.uni-heidelberg.de}

\author[J.~Kim]{Jaehoon Kim}
\address[J.~Kim]{Department of Mathematical Sciences, Korea Advanced Institute of Science and Technology, South Korea}
\email{jaehoon.kim@kaist.ac.kr}

\author[D.~K\"uhn]{Daniela K\"uhn}
\address[D.~K\"uhn]{School of Mathematics, University of Birmingham, United Kingdom}
\email{d.kuhn@bham.ac.uk}

\author[D.~Osthus]{Deryk Osthus}
\address[D.~Osthus]{School of Mathematics, University of Birmingham, United Kingdom}
\email{d.osthus@bham.ac.uk}

\thanks{The research leading to these results was partially supported by the European Research Council under the European Union's Seventh Framework Programme (FP/2007--2013) / ERC Grant 306349 (S.~Glock, J.~Kim and D.~Osthus), and the DFG, grant no.~JO 1457/1-1 (F.~Joos). The research was also partially supported by the EPSRC, grant no. EP/N019504/1, and by the Royal Society and the Wolfson Foundation (D.~K\"uhn).}

\date{\today}

\begin{document}

\begin{abstract}
The Oberwolfach problem, posed by Ringel in 1967, asks for a decomposition of $K_{2n+1}$ into edge-disjoint copies of a given $2$-factor.
We show that this can be achieved for all large $n$.
We actually prove a significantly more general result, which allows for decompositions into more general types of factors.
In particular, this also resolves the Hamilton-Waterloo problem for large $n$.
\end{abstract}
\maketitle

\section{Introduction}
A central theme in Combinatorics and related areas is the decomposition of large discrete objects into simpler or smaller objects.
In graph theory,
this can be traced back to the 18th century, when 
Euler asked for which orders orthogonal Latin squares exist (which was finally answered by Bose, Shrikhande, and Parker~\cite{BSP:60}).
This question can be reformulated as the existence question for resolvable triangle decompositions in the balanced complete tripartite graph.\COMMENT{A resolvable triangle decomposition of $K_{n,n,n}$ corresponds to a decomposition of a Latin square into transversals (rows, columns, entries).
Given a Latin square and a decomposition into transversals, we add the number of the transversal as a second coordinate to every entry and obtain a orthogonal Latin square.}
(Here a resolvable triangle decomposition is a decomposition into edge-disjoint triangle factors.)
In the 19th century,
Walecki proved the existence of decompositions of the complete graph $K_n$ (with $n$~odd) into edge-disjoint Hamilton cycles
and Kirkman formulated the school girl problem. The latter triggered the question for which $n$ the complete graph on $n$~vertices admits a resolvable triangle decomposition, which was finally resolved in the 1970s by Ray-Chaudhuri and Wilson~\cite{RW:71} and independently by Lu~\cite{lu:84} (see also~\cite{wilson:03}).
This topic has developed into a vast area with connections e.g.~to statistical design and scheduling, Latin squares and arrays,
graph labellings as well as combinatorial probability.

A far reaching generalisation of Walecki's theorem and Kirkman's school girl problem is the following problem posed by Ringel in Oberwolfach in~1967 (cf.~\cite{guy:71}).

\begin{problem}[Oberwolfach problem]
Let $n\in \N$ and let $F$ be a $2$-regular graph on $n$ vertices.
For which (odd) $n$ and $F$ does $K_{n}$ decompose into edge-disjoint copies of~$F$?
\end{problem}
 Addressing conference participants in Oberwolfach, Ringel fittingly formulated his problem as a scheduling assignment for diners: assume $n$ people are to be seated around round tables for $\frac{n-1}{2}$~meals, where the total number of seats is equal to $n$, but the tables may have different sizes. Is it possible to find a seating chart such that every person sits next to any other person exactly once?

We answer this affirmatively for all sufficiently large~$n$. We make no attempt to estimate or optimize the smallest $n$ for which our proof works. This makes our argument significantly clearer. It is known that for $n\le 60$, all but four instances have a solution~\cite{DFHMR:10,SDTBDC:19}.\COMMENT{$OP(3,3),OP(4,5), OP(3,3,5) ,OP(3,3,3,3)$. So for $n=6,9,11,12$}

A generalisation of the Oberwolfach problem is the Hamilton-Waterloo problem (cf.~\cite{ABBEZ:02}); here, two cycle factors are given and it is prescribed how often each of them is to be used in the decomposition. Informally, this means the conference takes place in two nearby venues (Hamilton and Waterloo), with possibly different seating arrangements.
We also resolve this problem in the affirmative (for large $n$) via the following even more general result. We allow an arbitrary collection of types of cycle factors, as long as one type appears linearly many times. This immediately implies that the Hamilton-Waterloo problem has a solution for large $n$ for any bounded number of given cycle factors.

\begin{theorem}\label{thm:multicity}
For every $\alpha>0$, there exists an $n_0\in \N$ such that for all odd $n\geq n_0$ the following holds.
Let $F_1,\ldots,F_k$ be $2$-regular graphs on $n$~vertices and let $m_1,\ldots,m_k\in \N$ be such that $\sum_{i\in [k]}m_i=(n-1)/2$ and $m_1\ge \alpha n$.
Then $K_n$ admits a decomposition into graphs $H_1,\ldots, H_{{(n-1)}/{2}}$ such that for exactly $m_i$ integers~$j$, the graph $H_j$ is isomorphic to~$F_i$.
\end{theorem}
Here we say a graph $G$ admits a \emph{decomposition} into graphs $H_1,\ldots, H_t$ if there exist edge-disjoint copies of $H_1,\ldots,H_t$ in $G$
such that every edge of $G$ belongs to exactly one copy.

Several authors (see e.g.~Huang, Kotzig, and Rosa~\cite{HKR:79}) considered a variant of the Oberwolfach problem for even~$n$; to be precise,
here we ask for a decomposition of $K_{n}$ minus a perfect matching
into $n/2-1$ copies of some given $n$-vertex $2$-regular graph~$F$.
We will deduce Theorem~\ref{thm:multicity} from a more general result (Theorem~\ref{thm:main}) which also covers this case.

The Oberwolfach problem and its variants have attracted the attention of many researchers, resulting in more than $100$~research papers covering a large number of partial results.
Most notably, Bryant and Scharaschkin~\cite{BS:09} proved it for infinitely many~$n$.
Traetta~\cite{traetta:13} solved the case when $F$ consists of two cycles only,
Bryant and Danziger~\cite{BD:11} solved the variant for even $n$ if all cycles are of even length, Alspach, Schellenberg, Stinson, and Wagner~\cite{ASSW:89} solved the case when all cycles have equal length (see~\cite{HS:91} for the analogous result for $n$ even), and Hilton and Johnson~\cite{HJ:01} solved the case when all but one cycle have equal length. The bipartite analogue of the Oberwolfach problem ($2$-factorizations of $K_{n,n}$) was solved in~\cite{piotrowski:91}. We refer to the surveys~\cite{alspach:96,BR:06} for more results on the Oberwolfach problem.

A related conjecture of Alspach from 1981 stated that for all odd $n$ the complete graph $K_n$ can be decomposed into any collection of cycles of length at most $n$ whose lengths sum up to~$\binom{n}{2}$.
This was solved by Bryant, Horsley, and Pettersson~\cite{BHP:14}.

Most classical results in the area are based on algebraic approaches,
often by exploiting symmetries.
More recently,  major progress for decomposition problems
has been achieved 
via absorbing techniques in combination with approximate decomposition results (often also in conjunction with probabilistic ideas).
This started off with
decompositions into Hamilton cycles~\cite{CKLOT:16,KO:13}, followed by the existence of combinatorial designs~\cite{GKLO:ta,keevash:14,keevash:18b} and progress on the tree packing conjecture~\cite{JKKO:19}.
In this paper, at a very high level,
we also pursue such an approach.
As approximate decomposition results, 
we exploit a hypergraph matching argument due to Alon and Yuster~\cite{AY:05} (which in turn is based on the R\"odl nibble via the Pippenger-Spencer theorem~\cite{PS:89})
and a bandwidth theorem for approximate decompositions due to Condon, Kim, K\"uhn, and Osthus~\cite{CKKO:19}.
Our absorption procedure utilizes as a key element a very special case of a recent result of Keevash on resolvable designs~\cite{keevash:18b}.

Whenever we only seek an approximate decomposition of a graph $G$,
the target graphs can be significantly more general and divisibility conditions disappear.
In particular, Allen, B{\"o}ttcher, Hladk{\'y}, and Piguet~\cite{ABHP:19} considered approximate decompositions into graphs of bounded degeneracy and maximum degree $o(n/\log n)$
whenever the host graph~$G$ is sufficiently quasirandom,
and Kim, K\"uhn, Osthus, and Tyomkyn~\cite{KKOT:19} considered 
approximate decompositions into graphs of bounded degree in host graphs~$G$ satisfying weaker quasirandom properties (namely, $\eps$-superregularity, see Section~\ref{sec:regularity}). 
Their resulting blow-up lemma for approximate decompositions was a key ingredient for~\cite{CKKO:19,JKKO:19} (and thus for the current paper too). It also implies that an approximate solution to the Oberwolfach problem can always be found (this was obtained independently by Ferber, Lee, and Mousset~\cite{FLM:17}).

Our Theorem~\ref{thm:multicity} actually follows from the following more general Theorem~\ref{thm:main}, which allows separable graphs.
An $n$-vertex graph $H$ is said to be \emph{$\xi$-separable} if there exists a set $S$ of at most $\xi n$ vertices such that every component of
$H \sm S$ has size at most~$\xi n$.
Examples of separable graphs include cycles, powers of cycles, planar graphs, and $F$-factors. More generally, for bounded degree graphs, the notion of separability
is equivalent to that of small bandwidth.

\begin{theorem} \label{thm:main}
For given $\Delta\in \mathbb{N}$ and $\alpha >0$, there exist $\xi_0>0$
and $n_0 \in \mathbb{N}$ such that the following holds for all $n\geq n_0$ and $\xi <\xi_0$.
Let $\cF, \cH$ be collections of graphs satisfying the following:
\begin{itemize}
\item $\cF$ is a collection of at least $\alpha n$ copies of $F$, where $F$ is a $2$-regular $n$-vertex graph;
\item each $H\in \cH$ is a $\xi$-separable  $n$-vertex $r_H$-regular graph for some $r_H\leq \Delta$;
\item $e(\cF\cup \cH)=\binom{n}{2}$.\COMMENT{equivalently: $2|\cF| + \sum_{H\in \cH}r_H = n-1$ }
\end{itemize}
Then $K_n$ decomposes into $\cF\cup \cH$.
\end{theorem}

Clearly, Theorem~\ref{thm:main} implies Theorem~\ref{thm:multicity} and also its corresponding version if $n$ is even and we ask for a decomposition of $K_n$ minus a perfect matching.

While far more general than the Oberwolfach problem,
Theorem~\ref{thm:main} may be just the tip of the iceberg, and it seems possible that
the following is true.
\begin{conj}
For all $\Delta \in \mathbb{N}$, there exists an $n_0 \in \mathbb{N}$ so that the following holds for all $n \ge n_0$.
Let $F_1,\ldots,F_t$ be $n$-vertex graphs such that $F_i$ is $r_i$-regular 
for some $r_i \le \Delta$ and $\sum_{i\in [t]}r_i={n-1}$.
Then there is a decomposition of $K_n$ into $F_1,\ldots,F_t$.
\end{conj}
The above conjecture is implicit in the `meta-conjecture' on decompositions proposed in~\cite{ABHP:19}.

Rather than considering decompositions of the complete graph $K_n$,
it is also natural to consider decompositions of host graphs of large minimum degree (this has applications e.g.~to the completion of partial 
decompositions of~$K_n$). 
Indeed, a famous conjecture of Nash-Williams~\cite{nash-williams:70} states that every $n$-vertex graph $G$
of minimum degree at least $3n/4$ has a triangle decomposition
(subject to the necessary divisibility conditions).
The following conjecture would (asymptotically) transfer this to arbitrary $2$-regular spanning graphs.

\begin{conj} \label{conj:mindegree}
Suppose $G$ is an $n$-vertex $r$-regular graph with even $r\geq \frac{3}{4}n+o(n)$
and $F$ is a $2$-regular graph on $n$ vertices.
Then $G$ decomposes into copies of~$F$.
\end{conj}
The (asymptotic version of the) Nash-Williams conjecture was reduced to its 
fractional version in~\cite{BKLO:16}. In combination with~\cite{dross:16}, this shows
that the Nash-Williams conjecture holds with $3n/4$ replaced by $9n/10+o(n)$.
There has also been considerable progress on decomposition problems involving
such host graphs of large minimum degree into other fixed subgraphs $H$ rather than triangles~\cite{GKLMO:19,montgomery:17}. It turns out that the chromatic number of $H$ is a crucial parameter for this problem. In particular, as proved in~\cite{GKLMO:19}, for bipartite graphs $H$ the `decomposition threshold' is always at most $\frac{2}{3}n +o(n)$.

Clearly, one can generalise Conjecture~\ref{conj:mindegree}
in this direction, e.g. to determine the decomposition threshold for $K_r$-factors. It might also be true that the `$3/4$' in Conjecture~\ref{conj:mindegree} can be replaced by `$2/3$' if $F$ consists only of even cycles.
We are confident that the ideas from this paper will be helpful in approaching these and other related problems. As remarked in Section~\ref{sec:the end}, the current method already allows us to replace $K_n$ by any $2r$-regular host graph $G$ with $2r\ge (1-\eps)n$.

This paper is organised as follows.
In Section~\ref{sec:proof sketch} we provide an overview of our approach
and in Section~\ref{sec:preliminaries} we collect several embedding and decomposition results including the bandwidth theorem for approximate decompositions from~\cite{CKKO:19} and a special case of a result of Keevash~\cite{keevash:18b} on resolvable designs. We prove Theorem~\ref{thm:main} in Section~\ref{sec:main} and add some concluding remarks in Section~\ref{sec:the end}.

\section{Proof sketch}  \label{sec:proof sketch}

For simplicity,
we just sketch the argument for the setting of the Oberwolfach problem;
that is, we aim to decompose $K_n$ into $\frac{n-1}{2}$ copies of an $n$-vertex $2$-regular graph $F$.
The proof essentially splits into two cases. 
In the first case we assume that almost all vertices of $F$ belong to `short' cycles, of length at most~$500$. Note that there must be some cycle length, say~$\ell^\ast$, such that at least $n/600$ vertices of $F$ lie in cycles of length~$\ell^\ast$. We will take a suitable number of random slices of the edges of~$K_n$ and then first embed, for every desired copy of $F$, all cycles whose lengths are different from~$\ell^\ast$. For this, we use standard tools based on the R\"odl nibble. We then complete the decomposition by embedding all the cycles of length~$\ell^\ast$. This last step uses a special case of a recent result of Keevash on the existence of resolvable designs~\cite{keevash:18b}.

The second case is much more involved and forms the core of the proof. We are now guaranteed that a (small) proportion of vertices of $F$ lies in `long' cycles.
To motivate our approach, consider the following simplified setup. Suppose $F$ consists only of cycles whose lengths are divisible by~$3$, and suppose for the moment we seek an $F$-decomposition of a $3$-partite graph $G$ with equitable vertex partition $(V_1,V_2,V_3)$ (so $G$ is a $C_3$-blowup). Let $\ell_1,\dots,\ell_t$ be the sequence of cycle lengths appearing in~$F$. Now, take any permutation $\pi$ on~$V_3$ which consists of cycles of lengths $\ell_1/3,\dots,\ell_t/3$. For instance, a $C_3$ in $F$ corresponds to a fixed point in~$\pi$, and a $C_6$ in $F$ corresponds to a transposition in~$\pi$. Now, define an auxiliary graph $\pi(G)$ by `rewiring' the edges between $V_2$ and $V_3$ according to~$\pi$. More precisely, we ensure that $E_{\pi(G)}(V_2,V_3)=\set{v_2\pi(v_3)}{v_2v_3\in E(G)}$. Suppose that $F'$ is a $C_3$-factor in $\pi(G)$. By `reversing' the rewiring, we obtain a copy of $F$ in~$G$. More precisely, let $\pi^{-1}(F')$ be the graph obtained from $F'$ by replacing $F'[V_2,V_3]$ with $\set{v_2\pi^{-1}(v_3)}{v_2v_3\in E(F')}$. Clearly, $\pi^{-1}(F')\cong F$ and $\pi^{-1}(F') \In G$. What is more, this rewiring is canonical in the following sense: if $F'$ and $F''$ are edge-disjoint $C_3$-factors in $\pi(G)$, then $\pi^{-1}(F')$ and $\pi^{-1}(F'')$ will be edge-disjoint copies of $F$ in~$G$. Thus, a resolvable $C_3$-decomposition of $\pi(G)$ immediately translates into an $F$-decomposition of~$G$.

Similarly, if all cycle lengths in $F$ are divisible by $4$, we can reduce the problem of finding an $F$-decomposition of a $C_4$-blowup to the problem of finding a resolvable $C_4$-decomposition of a suitably rewired $C_4$-blowup.
In order to deal with arbitrary $2$-regular graphs~$F$, we interweave such constructions for $C_3$, $C_4$ and~$C_5$. 
In Sections~\ref{subsec:absorbing gadget} and~\ref{subsec:rob dec}, we will construct an `absorbing graph' $G$ which is a partite graph on $18$~vertex classes such that finding an $F$-decomposition of~$G$ can be reduced to finding resolvable $C_3,C_4,C_5$-decompositions of suitable auxiliary graphs, in a similar way as sketched above. Crucially, $G$ has this property in a robust sense: even if we delete an arbitrary sparse graph $L$ from~$G$, as long as some necessary divisibility conditions hold, we are still able to find an $F$-decomposition of~$G-L$. 

The overall strategy is thus as follows: first, we remove $G$ from~$K_n$. Then we find an approximate decomposition of the remainder, which leaves a sparse leftover. 
For this, we employ the recent bandwidth theorem for approximate decompositions~\cite{CKKO:19}. 
(The existence of an approximate decomposition of $K_n-G$ would also follow directly from the blow-up lemma for approximate decompositions~\cite{KKOT:19},
but this would leave a leftover whose density is larger than that of the absorbing graph~$G$,
making our approach infeasible.)
We then deal with this leftover by using some edges of~$G$, in a very careful way, such that the remainder of $G$ is still appropriately divisible. The remainder of $G$ then decomposes as sketched above. In order to decompose the auxiliary graphs, we again use a very special case of the main result in~\cite{keevash:18b}.
The fact that we are guaranteed that $F$ has some long cycles will be helpful to construct the absorbing graph~$G$, more precisely, to ensure that all the $18$~vertex classes are of linear size. 
It is also essential when dealing with the leftover of the approximate decomposition.

\section{Preliminaries}\label{sec:preliminaries}

In this section,
we first introduce some notation,
then a Chernoff-type concentration inequality,
several graph embedding tools and notation concerning quasirandomness,
as well as the framework of the result of Keevash on resolvable decompositions.

\subsection{Notation}

For a $2$-regular graph $F$, let $\cC(F)$ denote the collection of cycles in~$F$. We also refer to a $2$-regular graph as a \defn{cycle factor}.
Let $G$ be a graph.
We denote by $e(G)$ the number of edges of~$G$, and by $|G|$ the number of vertices of~$G$.
For sets $U,U'\sub V(G)$,
we define $e_G(U)$ as the number of edges of the graph induced by $U$ and
$e_G(U,U')$ as the number of pairs $(u,u')\in U\times U'$ such that $uu'\in E(G)$. 
Hence, $e_G(U,U)=2e_G(U)$.
For a vertex $v\in V(G)$, we define $d_G(v,U):=|N_G(v)\cap U|$, where $N_G(v)$ is the neighbourhood of $v$ in~$G$. For a subgraph $H\In G$, we write $G-H$ for the graph with vertex set $V(G)$ and edge set $E(G)\sm E(H)$. We write $G\sm X$ for the graph obtained from $G$ by removing the vertices of~$X$, and $G\sm H:=G\sm V(H)$.
Given two graphs $G_1$ and $G_2$, define $G_1\bigtriangleup G_2$ to be the graph on $V(G_1)\cup V(G_2)$ whose edge set is $(E(G_1)\setminus E(G_2))\cup (E(G_2)\setminus E(G_1))$.

Given graphs $F$ and~$G$, a function $\sigma\colon V(F)\to V(G)$ is a \defn{homomorphism} if $\sigma(x)\sigma(y)\in E(G)$ for all $xy\in E(F)$. An injective homomorphism is called an \defn{embedding}. (Note that non-edges need not be preserved, that is, the corresponding subgraphs of $F$ in $G$ are not required to be induced.)
For a fixed graph $F$, an \defn{$F$-factor} in a graph $G$ is a collection $\cF$ of vertex-disjoint copies of $F$ in $G$ which cover all vertices of~$G$.
An \defn{$F$-decomposition} of $G$ is a collection $\cF$ of edge-disjoint copies of $F$ in $G$ which cover all edges of~$G$. An $F$-decomposition $\cF$ is called \defn{resolvable} if it can be partitioned into $F$-factors.

For a collection of graphs $\cH$, we define $e(\cH)$ by $\sum_{H\in \cH}e(H)$.
We write $G-\cH$ for the graph with vertex set $V(G)$ and edge set $E(G)\setminus \bigcup_{H\in \cH}E(H)$.
We also write $\Delta(\cH)$ for the maximum degree of $\bigcup_{H\in \cH}H$.
We say that $\cH=\Set{H_1,\dots,H_t}$ \defn{packs} into a graph $G$ if there exist edge-disjoint subgraphs $H_1',\dots,H_t'$ in $G$ such that $H_i'$ is a copy of $H_i$ for each $i\in[t]$.

Let $\cX$ be a set of disjoint vertex sets and $R$ a graph on~$\cX$. If $G$ is a graph with vertex partition~$\cX$, where each $X\in \cX$ is independent in~$G$, and such that $e_G(X,X')=0$ for all distinct $X,X'\in \cX$ with $XX'\notin E(R)$, then we say that $G$ has \defn{reduced graph}~$R$.

For a digraph $D$ and a vertex $v\in V(D)$,
we write $d^+_D(v)$ and $d^{-}_D(v)$ for the number of outgoing and incoming arcs at~$v$, respectively. We say that $D$ is $r$-regular if $d^+_D(v)=d^-_D(v)=r$ for all $v\in V(D)$.
Sometimes we write $\vv{D}$ for an oriented graph
and then $D$ is the undirected graph obtained from $\vv{D}$ by ignoring the orientations of $\vv{D}$.

We write $\N_0:=\N\cup \Set{0}$ and $[n]:=\Set{1,\dots,n}$. For a set of objects indexed by $[t]$, we often treat indices modulo~$t$, and define $a \mod{t}$ to be the unique integer $b\in [t]$ such that $a\equiv b \mod{t}$.
For $a,b,c\in \mathbb{R}$,
we write $a=b\pm c$ whenever $a\in [b-c,b+c]$.
For $a,b,c\in (0,1]$,
we write $a\ll b \ll c$ in our statements to mean that there are increasing functions $f,g:(0,1]\to (0,1]$
such that whenever $a\leq f(b)$ and $b \leq g(c)$,
then the subsequent result holds.

\subsection{Probabilistic tools}

At several stages of our proof we will apply the following standard Chernoff-type concentration inequalities.

\begin{lemma}[see {\cite[{Corollary~2.3, Corollary~2.4 and Theorem 2.8}]{JLR:00}}] \label{lem:chernoff}
Let $X$ be the sum of $n$ independent Bernoulli random variables. Then the following hold.
\begin{enumerate}[label={\rm(\roman*)}]
\item For all $0\le\eps \le 3/2$, we have $\prob{|X - \expn{X}| \geq \eps\expn{X} } \leq 2\eul^{-\eps^2\expn{X}/3}$.\label{chernoff eps}
\item If $t\ge 7 \expn{X}$, then $\prob{X\ge t}\le \eul^{-t}$.\label{chernoff crude}
\end{enumerate}
\end{lemma}

The following follows easily from Lemma~\ref{lem:chernoff}\ref{chernoff eps}.
An explicit derivation can be found in~\cite{GKLO:ta}.

\begin{lemma}\label{lem:separable chernoff}
Let $1/n\ll p,\alpha,1/B$. Let $\cI$ be a set of size at least $\alpha n$ and let $(X_i)_{i\in \cI}$ be a family of Bernoulli random variables with $\prob{X_i=1}\ge p$. Suppose that $\cI$ can be partitioned into at most $B$ sets $\cI_1,\dots,\cI_k$ such that for each $j\in[k]$, the variables $(X_i)_{i\in\cI_j}$ are independent. Let $X:=\sum_{i\in \cI}X_i$. Then we have
\begin{align*}
\prob{X \neq (1\pm n^{-1/5}) \expn{X}} \le \eul^{-n^{1/6}}.
\end{align*}
\end{lemma}

\subsection{Embedding and decomposition results}\label{sec:regularity}
Frequently in our proof we want to embed parts of a $2$-regular graph into `random-like' graphs.
For such a task, the blow-up lemma developed by Koml\'os, S\'ark\"ozy and Szemer\'edi~\cite{KSS:97} is a standard tool. Roughly speaking, it says that given a $k$-partite graph $G$ that is `super-regular' between any two vertex classes, and a $k$-partite bounded-degree graph $H$ with a matching vertex partition, then $H$ is a subgraph of~$G$. The notion `super-regular' is tailored towards being used after an application of Szemer\'edi's regularity lemma. Since we do not use Szemer\'edi's regularity lemma, but work essentially in random subgraphs of the complete graph, we can use a more convenient notion which is defined as follows.

We say that a graph $G$ on $n$ vertices is \defn{$(\eps,d)$-quasirandom} if $d_G(v)=(d\pm \eps)n$ for all $v\in V(G)$ and $|N_G(v_1)\cap N_G(v_2)|=(d^2\pm \eps)n$ for all distinct $v_1,v_2\in V(G)$.

If $V_1$ and $V_2$ are disjoint vertex sets in $G$, we also say that $G[V_1,V_2]$ is \defn{$(\eps,d)$-quasirandom} if for both $j\in[2]$, we have that $d_G(v,V_{3-j})=(d\pm \eps)|V_{3-j}|$ for all $v\in V_j$ and $|N_G(v_1)\cap N_G(v_2) \cap V_{3-j}|=(d^2\pm \eps)|V_{3-j}|$ for all distinct $v_1,v_2\in V_j$.
 It is well known that these conditions imply super-regularity (see~\cite{DLR:95}\COMMENT{Proposition 2.5}), that is, in addition to the degree condition, one also knows that between any large enough sets $V_1'\In V_1,V_2'\In V_2$, the edge density is very close to~$d$.
We need the following version of the blow-up lemma of Koml\'os, S\'ark\"ozy and Szemer\'edi, which we just state in the setting of quasirandom pairs.

\begin{lemma}[Blow-up lemma, \cite{KSS:97}] \label{lem:blow-up new}
Let $1/n \ll \eps \ll \kappa,d,1/\Delta,1/r$. Suppose that $G$ is an $n$-vertex graph with vertex partition $\cX$ and reduced graph $R$, where $|R|\le r$. Assume that $G[X,X']$ is $(\eps,d_{XX'})$-quasirandom for some $d_{XX'}\ge d$ whenever $XX'\in E(R)$. Assume also that $\min_{X\in \cX}|X|\ge \kappa \max_{X\in \cX}|X|$.

Let $H$ be a graph with $\Delta(H)\le \Delta$. Let $X_0\In V(H)$ be independent such that $|X_0|\le \eps n$ and no two vertices in $X_0$ have a common neighbour in~$H$.
Assume that $\sigma\colon H \to R$ is a homomorphism such that $|\sigma^{-1}(X)|=|X| $ for all $X\in \cX$, and $\phi_0\colon X_0 \to V(G)$ is an injective function such that $\phi_0(x)\in \sigma(x)$ for all $x\in X_0$.

Then there exists an embedding $\phi \colon H\to G$ which extends $\phi_0$ such that $\phi(x) \in \sigma(x)$ for all $x\in V(H)$.
\end{lemma}

The following is a straightforward consequence of the blow-up lemma.\COMMENT{If vertices in $X_0$ have common neighbours, embed them manually first and call the new fixed vertex set $X_0'$. Then apply blow-up lemma with $X_0'$. (Use Hajnal-Szemeredi or some other tool to partition $H$ into fairly equal-sized clusters, partition $G$ randomly with appropriate cluster sizes.)}

\begin{cor} \label{cor:quasirandom blowup}
Suppose $1/n\ll \eps \ll d$.
Let $G$ be an $(\eps,d)$-quasirandom $n$-vertex graph and suppose that $H$ is a graph on (at most) $n$~vertices with $\Delta(H)\le 2$. 
Suppose an independent set $X_0\In V(H)$ such that $|X_0|\leq \eps n$ and an injective function $\phi_0\colon X_0 \to V(G) $ are given.
Then there exists an embedding $\phi$ of $H$ into $G$ which extends~$\phi_0$.
\end{cor}

The following result due to Condon, Kim, K\"uhn and Osthus~\cite{CKKO:19} is a key ingredient in our approach.
It is in turn based on the blow-up lemma for approximate decompositions~\cite{KKOT:19}, and uses Szemer\'edi's regularity lemma.

\begin{theorem}[\cite{CKKO:19}]\label{thm: Padraig}
For all $\Delta\in \mathbb{N}\backslash\{1\}$, $0<\nu<1$, there exist $\xi >0$ and $n_0 \in \mathbb{N}$ such that for all $n \ge n_0$ and $1 -\frac{1}{200\Delta} \le d \leq 1$ the following holds:
Suppose that $\mathcal{H}$ is a collection of $n$-vertex $\xi$-separable graphs and $G$ is an $n$-vertex graph such that
\begin{enumerate}[label=\rm{(\roman*)}]
\item $d_G(x)=(d \pm \xi) n$ for all $x\in V(G)$;
\item $\Delta(H) \le \Delta$ for all $H \in \mathcal{H}$;
\item $e(\mathcal{H}) \le (1-\nu)e(G)$.
\end{enumerate}
Then $\mathcal{H}$ packs into $G$.
\end{theorem}

\subsection{Partite decompositions of typical graphs}

We will make use of a recent result of Keevash~\cite{keevash:18b} on partite decompositions of typical (hyper-)graphs. His result applies in a far more general setting and we only need a simple consequence thereof here (see Theorem~\ref{thm:keevash partite}). For simplicity, we only introduce the relevant concepts for partite graph decompositions.

We first define typicality. Roughly speaking, a graph is typical if common neighbourhoods are as large as one would expect in a random graph. We also need such a notion for partite structures. Suppose that $G$ is a graph with vertex partition $(V_1,\dots,V_t)$ (we do not assume that these sets are independent). Let $\tau\colon V(G)\to [t]$ denote the assignment function such that $v\in V_{\tau(v)}$ for all $v\in V(G)$. Given a symmetric matrix $D\in [0,1]^{t\times t}$, $s\in \N$ and $\eps>0$, we say that $G$ is \defn{$(\eps,s,D)$-typical} if for any set $S\In V(G)$ with $|S|\le s$ and any $i \in [t]$, we have that
\begin{align}
\left|V_i \cap \bigcap_{v\in S}N_G(v)\right| &= (1\pm \eps)  |V_i| \prod_{v\in S} D_{\tau(v)i}. \label{eq:typicality def}
\end{align}
Note that if $G$ is a weighted  binomial random graph where an edge between $V_i$ and $V_j$ is included with probability $D_{ij}$, then \eqref{eq:typicality def} holds with high probability. Note also that if $D_{ij}=0$ and \eqref{eq:typicality def} holds, then $e_G(V_i,V_j)=0$. 
We will often write $D_{V_iV_j}$ instead of $D_{ij}$ and $D_{V_i}$ instead of $D_{V_iV_i}$.
If $t=1$, we simply write $(\eps,s,D_{11})$-typical instead of $(\eps,s,D)$-typical.

Our aim is to decompose $G$ into a given graph~$H$, following a prescribed pattern. More precisely, let $\sigma\colon V(H)\to [t]$ be an assignment of the vertices of $H$ to the vertex partition classes of~$G$. Given an embedding $\phi\colon H\to G$ of $H$ into $G$ such that $\tau(\phi(x))=\sigma(x)$ for all $x\in V(H)$, we say that $\phi(H)$ is a \defn{$\sigma$-copy of $H$ in~$G$.} A collection $\cH$ of edge-disjoint $\sigma$-copies of $H$ in $G$ is called an \defn{$(H,\sigma)$-packing in~$G$}. We say that $\cH$ is an \defn{$(H,\sigma)$-decomposition of~$G$} if, in addition, every edge of $G$ is covered.

Roughly speaking, Theorem~\ref{thm:keevash partite} guarantees an $(H,\sigma)$-decomposition of $G$ under two assumptions: typicality and divisibility. We have already defined typicality of~$G$. Now we relate the density matrix $D$ to the given assignment $\sigma$ of~$H$. Let $(P_1,\dots,P_t)$ be the partition of $V(H)$ induced by $\sigma$, i.e.~$P_i=\sigma^{-1}(i)$ for each $i\in[t]$. We let $I^{H,\sigma}$ denote the symmetric indicator $(t\times t)$-matrix defined as $I^{H,\sigma}_{ij}:=1$ if $e_H(P_i,P_j)>0$ and $I^{H,\sigma}_{ij}:=0$ otherwise.
We will later require that $D\ge d\cdot I^{H,\sigma}$. This accounts for the fact that if some edge of $H$ is mapped to $(V_i,V_j)$, then we require this pair to be sufficiently dense in~$G$.

Finally, we say that $G$ is \defn{$(H,\sigma)$-divisible} if the following hold:
\begin{enumerate}[label=\rm{(\roman*)}]
\item there exists $m\in \N_0$ such that $e_G(V_i,V_j)=m \cdot e_H(P_i,P_j)$ for all $i,j\in [t]$;
\item for all $i\in [t]$ and every $v\in V_i$, there are $(a^v_x)_{x\in P_i}\in \N_0$ such that $d_G(v,V_j)=\sum_{x\in P_i}a^v_x \cdot d_H(x,P_j)$ for all $j\in[t]$.
\end{enumerate}
It is easy to see that $(H,\sigma)$-divisibility is a necessary condition for an $(H,\sigma)$-decomposition to exist.\COMMENT{$m$ is number of copies in decomposition, $a^v_x$ is the number of times $v$ plays the role of $x$ in a $\sigma$-copy of $H$}
If $\sigma$ is bijective, we follow the notation in \cite{keevash:18b} and simply say that $G$ is \defn{$H$-balanced} if it is $(H,\sigma)$-divisible, and an $(H,\sigma)$-decomposition is simply called a \defn{partite $H$-decomposition}.

\begin{theorem}[{\cite[cf.~Theorem~7.8]{keevash:18b}}]  \label{thm:keevash partite}
Let $1/n\ll \eps \ll 1/s \ll d, 1/h$. Let $G$ be $(\eps,s,D)$-typical with vertex partition $(V_1,\dots,V_t)$ such that $d n \le |V_i|\le n$ for all $i\in[t]$.
Let $H$ be a graph on $h$ vertices and $\sigma\colon V(H)\to [t]$. Suppose that $D\ge d\cdot I^{H,\sigma}$ and that $G$ is $(H,\sigma)$-divisible. Then $G$ has an $(H,\sigma)$-decomposition.
\end{theorem}

We now use Theorem~\ref{thm:keevash partite} to deduce the two results about resolvable cycle decompositions which we will need later.
The following will be used in Case~1 of our proof, where most vertices of $F$ are contained in cycles of length at most~$500$.

\begin{cor} \label{thm:non partite cycle dec}
Let $1/n\ll \eps \ll 1/s \ll d, 1/\ell$.
Let $G$ be $(\eps,s,D)$-typical with vertex partition $(V,U)$ such that $d n\le |V|,|U| \le n$. Assume that $D_{V},D_{VU}\ge d$ and $D_{U}=0$.
Let $W_\ell$ be the wheel graph with $\ell$ spokes and hub $w$, and let $\sigma$ assign $w$ to $U$ and all other vertices to~$V$.
Assume that $d_G(v,V)=2d_G(v,U)$ for all $v\in V$ and $\ell \mid d_G(u)$ for all $u\in U$.
Then $G$ has a $(W_\ell,\sigma)$-decomposition.
\end{cor}

\proof
Note that $G$ is $(W_\ell,\sigma)$-divisible,\COMMENT{we have $e_G(V)=e_G(V,U)$ and $e_H(V(W_\ell)-w)=e_H(V(W_\ell)-w,w)=\ell$ and hence $$\frac{e_G(V)}{\ell}=\frac{e_G(V,U)}{\ell}=\frac{\sum_{u\in U}d_G(u)}{\ell}$$ is an integer. The degree conditions clearly hold.} so we can apply Theorem~\ref{thm:keevash partite}.
\endproof

Note that given such a decomposition, for every vertex $u\in U$, the collection of all cycles which together with $u$ form a wheel in the decomposition form a $C_\ell$-factor of $G[N_G(u)]$.

We will also need the following approximate version of Corollary~\ref{thm:non partite cycle dec}, which is much simpler to prove and follows from standard hypergraph matching results based on the R\"odl nibble. Note that we do not need to assume divisibility in this case, which makes it more convenient to apply.

\begin{cor} \label{cor:approximate wheel}
Let $1/n\ll \eps \ll \gamma, d, 1/\ell$. Let $G$ be $(\eps,\ell,D)$-typical with vertex partition $(V,U)$ such that $d n\le |V|,|U| \le n$. Assume that $D_{V},D_{VU}\ge d$ with $D_{UV}|U|=(1\pm \eps)|D_V||V|/2$,\COMMENT{implies that $e_G(U,V) \approx e_G(V)$.} and $D_{U}=0$.
Let $W_\ell$ be the wheel graph with $\ell$ spokes and hub $w$, and let $\sigma$ assign $w$ to $U$ and all other vertices to~$V$.
Then $G$ has a $(W_\ell,\sigma)$-packing such that the leftover $L$ satisfies $\Delta(L)\le \gamma n$.
\end{cor}

\proof
Define an auxiliary $2\ell$-uniform hypergraph $\cH$ with vertex set $E(G)$ where the edges of $\cH$ correspond to $\sigma$-copies of $W_\ell$ in~$G$. Using the typicality condition, we can count that every edge of $G$ between $V$ and $U$ lies in $(1\pm \eps)^{\ell-1} \cdot \frac{1}{2} D_V^{\ell}D_{UV}^{\ell-1} |V|^{\ell-1}$ $\sigma$-copies of $W_\ell$,\COMMENT{factor $\frac{1}{2}$ because can go around the cycle in two ways} and every edge contained in $V$ lies in $(1\pm \eps)^{\ell-1} D_V^{\ell-1}D_{UV}^\ell |U||V|^{\ell-2}$ such copies. Using $D_{UV}|U|=(1\pm \eps)D_V |V|/2$, we can conclude that $d_\cH(e)= (1\pm \eps)^{\ell} \frac{1}{2} D_V^{\ell}D_{UV}^{\ell-1} |V|^{\ell-1}$ for all $e\in V(\cH)$, i.e.~$\cH$ is almost regular. For $v\in V(G)$, let $F_v\In V(\cH)$ be the set of edges of $G$ which are incident to~$v$.
A result of Alon and Yuster~\cite[Theorem~1.2]{AY:05} implies that there exists a matching $M$ in $\cH$ such that for each set $F_v$, all but at most $\gamma |F_v|/2$ vertices of $F_v$ are covered by~$M$. Clearly, $M$ corresponds to a $(W_\ell,\sigma)$-packing in $G$ such that the leftover graph $L$ satisfies $\Delta(L)\le \gamma n$.
\endproof

The next result will be used in Case~2 of our proof. It asserts the existence of resolvable partite cycle decompositions in typical partite graphs. In the proof, we add a new vertex class of size equal to the number of cycle factors required for a resolvable decomposition, and join it completely to the rest of the graph. A wheel decomposition of this auxiliary graph encodes a resolvable cycle decomposition of the original graph. A very similar reduction has also been used e.g.~in~\cite{keevash:18b}\COMMENT{similar ideas used by Lamken and Wilson} to derive the existence of resolvable designs. We include a short proof for completeness.

\begin{cor} \label{thm:res partite cycle dec}
Let $1/n\ll \eps \ll 1/s \ll d, 1/\ell$. Let $G$ be $(\eps,s,D)$-typical with vertex partition $(V_1,\dots,V_\ell)$ into sets of size $n$ each, where $D_{i,i+1}:=D_{i+1,i}:=d$ for all $i\in[\ell]$ (indices modulo $\ell$) and $D_{ij}=0$ otherwise.
Assume that there exists $r\in \N$ such that $d_G(v,V_{i-1})=d_G(v,V_{i+1})=r$ for all $v\in V_i$ and $i\in[\ell]$ (indices modulo $\ell$).\COMMENT{From typicality, it follows that $r=(1\pm \eps)dn$} Then $G$ has a resolvable partite $C_\ell$-decomposition.
\end{cor}

\proof
Let $W_\ell$ be the wheel graph with $\ell$ spokes and hub $w$.
We extend $G$ to a graph $G'$ by adding a new vertex class $V_w$ of size $r$ and joining each $v\in V_w$ to all vertices of~$G$. Accordingly, extend $D$ to a $((\ell+1)\times(\ell+1))$-matrix $D'$ by defining $D'_{V_w}:=0$ and $D'_{V_wV_i}:=D'_{V_iV_w}:=1$ for all $i\in[\ell]$. Noting that $V_i\In N_{G'}(v)$ for all $i\in[\ell]$ and $v\in V_w$, it is easy to see that the typicality of $G$ directly implies that $G'$ is $(\eps,s,D')$-typical.

It is also straightforward to check that $G'$ is $W_\ell$-balanced.\COMMENT{for a vertex in an original cluster, all degrees to neighbouring clusters are equal to $r$. For a new vertex, all degrees to the other clusters are $n$.}
Thus, applying Theorem~\ref{thm:keevash partite} to $G'$ yields a partite $W_\ell$-decomposition $\cW$ of~$G'$. This gives us a resolvable partite $C_\ell$-decomposition $\cC$ of $G$ as follows. For each vertex $v\in V_w$, let $\cW_v$ be the set of all copies of $W_\ell$ in $\cW$ which contain~$v$. Let $\cC_v$ be obtained from $\cW_v$ by removing $v$ from each element. Then $\cC_v$ is a $C_\ell$-factor of $G$,\COMMENT{the center of the wheel always lies in $V_w$} and $\cC:=\bigcup_{v\in V_w}\cC_v$ is a resolvable partite $C_\ell$-decomposition of~$G$.
\endproof

The following proposition states that random slices of typical graphs are again typical. This is an easy consequence of Lemma~\ref{lem:chernoff}\ref{chernoff eps}.

\begin{prop} \label{prop:typical random slice}
Let $1/n\ll \eps, 1/s,d,p$.
Let $G$ be $(\eps,s,D)$-typical with vertex partition $(V_1,\dots,V_t)$ such that $d n  \le |V_i|\le n$ for all $i\in[t]$ and $D \ge d I$ for some indicator $(t\times t)$-matrix~$I$. Suppose we choose a random subgraph $G'$ of $G$ by including each edge independently with probability~$p$. Then $G'$ is $(1.1\eps,s,pD)$-typical with probability at least $1-\eul^{-\sqrt{n}}$.
\end{prop}
\COMMENT{
\proof Consider a set $S$ of at most $s$ vertices of $G$ and some $i\in[t]$. We can clearly assume that $D_{\tau(v)i}\ge d$ for all $v\in S$. For every vertex $u\in V_i \cap \bigcap_{v\in S}N_G(v)$, let $X_u$ be the indicator variable of the event that $vu\in E(G')$ for all $v\in S$. Clearly, $\prob{X_u}=p^{|S|}$. Thus
$$\expn{\left|V_i \cap \bigcap_{v\in S}N_{G'}(v)\right|} = p^{|S|}\left|V_i \cap \bigcap_{v\in S}N_G(v)\right| = (1\pm \eps)  |V_i| \prod_{v\in S} pD_{\tau(v)i}.$$
Moreover, the Bernoulli variables $X_u$ are independent. Thus, Lemma~\ref{lem:chernoff}\ref{chernoff eps} implies that $$\prob{\left|V_i \cap \bigcap_{v\in S}N_{G'}(v)\right|\neq (1\pm 1.1\eps)  |V_i| \prod_{v\in S} pD_{\tau(v)i}} \le \eul^{-n^{0.9}}.$$
A union bound over all choices of $S$ and $i$ then yields the result.
\endproof
}

We conclude this subsection with the following simple fact about the robustness of the typicality property.

\begin{prop} \label{prop:typical noise}
Let $G$ be $(\eps,s,D)$-typical with vertex partition $(V_1,\dots,V_t)$ such that $d n  \le |V_i|\le n$ for all $i\in[t]$ and $D \ge d I$ for some indicator $(t\times t)$-matrix~$I$. Suppose that $L$ is a graph on $V(G)$ such that whenever $L[V_i,V_j]$ is non-empty, then neither is~$G[V_i,V_j]$. Suppose that $\Delta(L)\le \gamma n$. Then $G\bigtriangleup L$ is still $(\eps+s\gamma d^{-s-1},s,D)$-typical.
\end{prop}

\proof Consider a set $S$ of at most $s$ vertices of $G$ and some $i\in[t]$. If $D_{\tau(v)i}=0$ for some $v\in S$, then $d_{G\bigtriangleup L}(v,V_i)=0$ and hence there is nothing to show. So we can assume that $D_{\tau(v)i}\ge d$ for all $v\in S$.
We then have $|V_i \cap \bigcap_{v\in S}N_{G\bigtriangleup L}(v)|= |V_i \cap \bigcap_{v\in S}N_G(v)|\pm s\gamma n$ and $s\gamma n\le s\gamma  d^{-s-1} |V_i| \prod_{v\in S} D_{\tau(v)i}$.
\endproof

\section{Main proof}  \label{sec:main}

In this section, we prove Theorem~\ref{thm:main}. 
The proof divides into two cases. 
The first case assumes that the given $2$-regular graph $F$ has very few vertices in long cycles. In the second case, we are guaranteed that linearly many vertices of $F$ lie in long cycles. The latter case is much more involved than the first case and needs some additional preliminary work. In the following subsections, we will develop the necessary tools for Case~2.

The following definition will be used throughout the section.
Let $\cX$ be a set of disjoint vertex sets and $R$ an oriented graph on~$\cX$. Assume that $G$ is a graph with vertex partition~$\cX$ and reduced graph~$R$ (where each $X\in \cX$ is independent in~$G$). Let $G_{R}$ be the oriented graph obtained from $G$ by orienting every edge $e\in E(G)$ with the same direction as the reduced edge of~$R$ corresponding to~$e$.

We say that $G$ is \defn{$r$-balanced with respect to $(\cX,R)$} if $G_R$ is $r$-regular, that is, $d^+_{G_R}(v)=d^-_{G_R}(v)=r$ for all $v\in V(G)$.
We simply say that $G$ is \defn{balanced} if it is $r$-balanced for some $r\in \N$.

\subsection{Cyclic partitions}  \label{subsec:absorbing gadget}

As sketched in Section~\ref{sec:proof sketch}, we reduce the problem of finding an $F$-decomposition of a graph to finding resolvable $C_3,C_4,C_5$-decompositions of suitable auxiliary graphs. 

For $\ell\in \N$, we say that $(a_1,\dots,a_t)$ is a \defn{cyclic partition of $\ell$} if $a_i\in \N$ for all $i\in[t]$ and $\sum_{i\in [t]}a_i = \ell$. We identify $(a_1,\dots,a_t)$ with $(a_i,\dots,a_t,a_1,\dots,a_{i-1})$ for each $i\in[t]$ and treat indices modulo~$t$. Moreover, for $S\In \N$, we say that $(a_1,\dots,a_t)$ is a \defn{cyclic $S$-partition of $\ell$} if $a_i\in S$ for all $i\in[t]$.

For a cyclic partition $\mathfrak{a}=(a_1,\dots,a_t)$ and a sequence $a'=(a_1',\dots,a_{t'}')$, 
we let $\fc^{\mathfrak{a}}(a')$ denote the number of appearances of $a'$ in~$\mathfrak{a}$, that is, 
$$\fc^{\mathfrak{a}}(a'):=|\set{i\in [t]}{a_{i+j \mod{t}}=a'_{j} \mbox{ for all }j\in[t']}|.$$
Note that if e.g.~$\mathfrak{a}=(a)$ and $a'=(a,a)$, we have $\fc^{\mathfrak{a}}(a')=1$.
To improve readability, we write $\fc^\fa(a)$ instead of $\fc^\fa((a))$ and $\fc^\fa(a,b)$ instead of $\fc^\fa((a,b))$.
For a sequence~$a$, we let $a^m$ denote the sequence which is the concatenation of $m$ copies of~$a$.

In this paper, we will only consider cyclic $\Set{3,4,5}$-partitions $\fa$
where $\fc^\fa(a,b)=\fc^\fa(b,a)$ for all $a,b\in \{3,4,5\}$.
For brevity, we will simply call these the \defn{admissible} partitions. 
Let $$\cI:= \{3,4,5\} \times \{3,4,5\}.$$

We now make some easy observations regarding admissible partitions.

\begin{prop}\label{prop:number partition}
For each $\ell \in \mathbb{N}$ with $\ell\ge 3$, 
there exists an admissible partition of~$\ell$. 
Moreover, if $\ell \ge 500$, then there exists an admissible partition $\fa$ of $\ell$ such that $\fc^\fa((a,b)^6)\ge \ell/200$ for all $(a,b)\in \cI$.
\end{prop}

\proof
For the first part, write $\ell=3x+4y$ with $x,y\in \Z$. If $\ell \ge 12$, we can clearly assume $x,y\ge 0$, and this can also easily be checked if $\ell <12$ and $\ell\neq 5$.
Thus, unless $\ell=5$, we obtain a partition $\fa$ consisting of $x$ $3$'s and $y$ $4$'s (in any order). 
Since $\fc^\fa(3,4)=\fc^\fa(4,3)$, the partition $\fa$ is also admissible.
If $\ell=5$, we simply take $(5)$. (This is in fact the reason why we allow $5$'s in admissible partitions, otherwise $3$ and $4$ would suffice.)

Now, assume that $\ell\ge 500$. Choose $\ell - 61 < \ell' \le \ell-13$ such that $\ell'$ is divisible by $48$. Then there is a partition of $\ell'$ such that each number $3,4,5$ appears exactly $\ell'/12$ times. 
Moreover, similarly as in the first part, there is an admissible partition $\fa'$ of $\ell-\ell'\ge 13$
that consists of a positive number of $4$'s followed by a positive number of $3$'s.
We now construct $\fa$ by concatenating $a^{\ell'/48}$ for each $a\in \Set{(3,3),(3,4),(4,4),(4,5),(5,5),(3,5)}$ (in this order), and finally attaching~$\fa'$.
It is easy to check that $\fa$ is an admissible partition of~$\ell$. 
Moreover, by construction, we have $\fc^\fa((a,b)^{6})\ge \ell'/48 - 6   \ge \ell/200$ for all $(a,b)\in \cI$.
\endproof

\begin{fact} \label{fact:simple partition facts}
For any admissible partition $\fa$ of~$\ell\geq 3$ the following hold:
\begin{enumerate}[label=\rm{(\roman*)}]
\item $\sum_{a\in\Set{3,4,5}} a\cdot \fc^\fa(a)= \ell$;\label{prop:simple partition facts:sum}
\item for each $a\in\Set{3,4,5}$, we have $\fc^\fa(a)=\sum_{b\in \Set{3,4,5}} \fc^\fa(a,b)$. \label{prop:simple partition facts:handshake}

\end{enumerate}
\end{fact}

For every $\ell\in \N$ with $\ell\ge 3$, let $\mathfrak{a}^\ell$ be an admissible partition of~$\ell$. 
Moreover, if $\ell\ge 500$, we also assume that $\fc^\fa((a,b)^6)\ge \ell/200$ for all $(a,b)\in \cI$. 
We fix $\fa^\ell$ throughout the remainder of the paper. Moreover, whenever $C$ is a cycle of length~$\ell$, we let $\fa^C:=\mathfrak{a}^\ell$ and $\fc^C(\cdot):=\fc^{\mathfrak{a}^C}(\cdot)$.
Hence, given any $2$-regular graph~$F$, we obtain a collection $$\mathfrak{a}^F:=\set{\mathfrak{a}^C}{C \in \cC(F)}$$ of cyclic partitions. Since we fix our choice of $\mathfrak{a}^\ell$ throughout, this collection is unique for every $2$-regular graph~$F$.
We define the counting function
\begin{align}
	\fc^{F}(\cdot):=\sum_{C\in \cC(F)}\fc^{C}(\cdot). \label{appearance counting factor}
\end{align}

\begin{prop} \label{prop:gadget sizes}
Let $F$ be a $2$-regular graph on $n$ vertices. Then the following hold:
\begin{enumerate}[label=\rm{(\roman*)}]
\item $\fc^F(a,b) = \fc^F(b,a)$ for all $(a,b)\in \cI$;\label{prop:gadget sizes:symmetry}
\item $\sum_{a\in \Set{3,4,5}} a\cdot \fc^F(a)= n$; \label{prop:gadget sizes:sum}
\item for each $a\in\Set{3,4,5}$, we have $\fc^F(a)=\sum_{b\in \Set{3,4,5}} \fc^F(a,b)$; \label{prop:gadget sizes:subpartition}
\item if at least $\eta n$ vertices of $F$ lie in cycles of length at least~$500$, then $\fc^F((a,b)^6) \ge \eta n/200 $ for all $(a,b)\in \cI$.\label{prop:gadget sizes:rich} 
\end{enumerate}
\end{prop}

\proof
Items \ref{prop:gadget sizes:symmetry}, \ref{prop:gadget sizes:sum} and \ref{prop:gadget sizes:subpartition} follow directly from Fact~\ref{fact:simple partition facts} and \eqref{appearance counting factor}. 
 \ref{prop:gadget sizes:rich} follows by assumption on $\fa^\ell$ and~\eqref{appearance counting factor}.\COMMENT{Let $\cC^{\ge 500}$ be the collection of all cycles in $F$ of length at least~$500$. 
Then for every $(a,b)\in \cI$, we have 
$$\fc^F((a,b)^6) \ge \sum_{C\in \cC^{\ge 500}}\fc^{C}((a,b)^6) \ge  \sum_{C\in \cC^{\ge 500}}|C|/200 \ge \eta n/200, $$ as required.}
\endproof

\begin{figure}[t]
\centering

\tikzstyle{circ}=[ultra thick,fill=white]
\tikzstyle{circc}=[ultra thick,fill=white]
\tikzstyle{lineblack}=[line width=4,black,postaction={decorate}]
\tikzstyle{linehalfblack}=[line width=4,black!50,postaction={decorate}]

\begin{tikzpicture}[scale=0.5]

\def\radc{2.5}

\draw[lineblack,decoration={markings, mark=at position 0.6 with {\arrow{>}}}](0,-2.5) .. controls (1,-1) and (1,1) .. (0,2.5);
\draw[lineblack,decoration={markings, mark=at position 0.75 with {\arrow{>}}}](0,2.5) -- (-.5,0);
\draw[lineblack,decoration={markings, mark=at position 0.6 with {\arrow{>}}}](-10,2.5).. controls (-9,3.5) and (-6,3.5) .. (-4.5,0) ;
\draw[lineblack,decoration={markings, mark=at position 0.60 with {\arrow{>}}}] (-10,2.5)--(-6.5,0);
\draw[lineblack,decoration={markings, mark=at position 0.60 with {\arrow{>}}}] (11,2.5)--(6.5,0);
\draw[lineblack,decoration={markings, mark=at position 0.6 with {\arrow{>}}}](11,2.5) .. controls (7,3.5) and  (6,3.5).. (4.5,0);
\draw[lineblack,decoration={markings, mark=at position 0.6 with {\arrow{>}}}] (0,2.5)--(-2.5,0);
\draw[lineblack,decoration={markings, mark=at position 0.7 with {\arrow{>}}}] (0,2.5) -- (2.5,0);

\draw[linehalfblack,decoration={markings, mark=at position 0.6 with {\arrow{>}}}] (4.5,0) .. controls (3,-2.5) and (2,-2.5) .. (0,-2.5);
\draw[linehalfblack,decoration={markings, mark=at position 0.6 with {\arrow{>}}}] (2.5,0) .. controls (4.5,-2.5) and (6.5,-2.5) .. (11,-2.5);
\draw[linehalfblack,decoration={markings, mark=at position 0.7 with {\arrow{>}}}](6.5,0) .. controls  (5,-6) and (-8,-7) .. (-10,-2.5);
\draw[linehalfblack,decoration={markings, mark=at position 0.6 with {\arrow{>}}}](-4.5,0) .. controls  (-3,-2.5) and (-2,-2.5) .. (0,-2.5);
\draw[linehalfblack,decoration={markings, mark=at position 0.7 with {\arrow{>}}}] (-2.5,0) .. controls (-5,-2.5) and (-8.5,-4) .. (-10,-2.5);
\draw[linehalfblack,decoration={markings, mark=at position 0.7 with {\arrow{>}}}] (-6.5,0).. controls (-2,-6) and (6,-6) .. (11,-2.5);
\draw[linehalfblack,decoration={markings, mark=at position 0.58 with {\arrow{>}}}](-.5,0) -- (0,-2.5);

\draw[circ] (0,-2.5) circle (1);
\draw[circ] (0,2.5) circle (1);
\draw[circc] (-0.5,0) circle (0.7);

\begin{scope}[shift={(11,0)}]

\draw[lineblack,decoration={markings, mark=at position 0.65 with {\arrow{>}}}](0:\radc)--(90:\radc);
\draw[lineblack,decoration={markings, mark=at position 0.75 with {\arrow{>}}}](90:\radc)--(180:\radc);
\draw[lineblack,decoration={markings, mark=at position 0.65 with {\arrow{>}}}](270:\radc)--(0:\radc);
\draw[linehalfblack,decoration={markings, mark=at position 0.65 with {\arrow{>}}}](180:\radc)--(270:\radc);

\draw[circ] (0:\radc) circle (1);
\draw[circ] (90:\radc) circle (1);
\draw[circc] (180:\radc) circle (0.7);
\draw[circ] (270:\radc) circle (1);

\draw (1.4,-3.5) node {$X_2^4$};
\draw (1.4,3.5) node {$X_4^4$};
\draw (3,-1.5) node {$X_3^4$};

\end{scope}

\begin{scope}[shift={(-11,0)}]

\draw (-0.5,-3.5) node {$X_2^5$};
\draw (-0.5,3.5) node {$X_5^5$};
\draw (-2.3,-3.1) node {$X_3^5$};
\draw (-2.3,3.1) node {$X_4^5$};

\draw[lineblack,decoration={markings, mark=at position 0.65 with {\arrow{>}}}](288:\radc)--(216:\radc);
\draw[lineblack,decoration={markings, mark=at position 0.65 with {\arrow{>}}}](216:\radc)--(144:\radc);
\draw[lineblack,decoration={markings, mark=at position 0.65 with {\arrow{>}}}](144:\radc)--(72:\radc);
\draw[lineblack,decoration={markings, mark=at position 0.7 with {\arrow{>}}}](72:\radc)--(0:\radc);
\draw[linehalfblack,decoration={markings, mark=at position 0.65 with {\arrow{>}}}](0:\radc)--(288:\radc);

\draw[circc] (0:\radc) circle (0.7);
\draw[circ] (72:\radc) circle (1);
\draw[circ] (144:\radc) circle (1);
\draw[circ] (216:\radc) circle (1);	
\draw[circ] (288:\radc) circle (1);	

\end{scope}

\draw[circc] 
(2.5,0) circle (0.7)
(4.5,0) circle (0.7)
(6.5,0) circle (0.7)
(-2.5,0) circle (0.7)
(-4.5,0) circle (0.7)
(-6.5,0) circle (0.7)
;

\draw (-10.3,0) node {$\pi$};
\draw[->,thick] (-9.7,0.4) arc (155:205:1);

\fill[nearly opaque,pattern=grid,fill=blue!30!black] 
(-0.5,0) circle (0.7)
(-2.5,0) -- (-3.2,0) arc (180:0:0.7) -- cycle
(2.5,0) -- (1.8,0) arc (180:0:0.7) -- cycle
(-4.5,0) -- (-5.2,0) arc (180:360:0.7) -- cycle
(4.5,0) -- (3.8,0) arc (180:360:0.7) -- cycle;

\fill[nearly opaque,pattern=grid,fill=black!30!green] 
(8.5,0) circle (0.7)
(2.5,0) -- (1.8,0) arc (180:360:0.7) -- cycle
(4.5,0) -- (3.8,0) arc (180:0:0.7) -- cycle
(6.5,0) -- (7.2,0) arc (0:180:0.7) -- cycle
(-6.5,0) -- (-7.2,0) arc (180:360:0.7) -- cycle
;

\fill[nearly opaque,pattern=grid,fill=orange!80!white] 
(-8.5,0) circle (0.7)
(-2.5,0) -- (-1.8,0) arc (360:180:0.7) -- cycle
(-4.5,0) -- (-3.8,0) arc (0:180:0.7) -- cycle
(-6.5,0) -- (-7.2,0) arc (180:0:0.7) -- cycle
(6.5,0) -- (7.2,0) arc (360:180:0.7) -- cycle
;

\draw[fill] (4.5,0) circle (0.1);

\draw[dashed, very thick]
(-1.8,0)--(-3.2,0)
(-3.8,0)--(-5.2,0)
(-5.8,0)--(-7.2,0)
(-7.8,0)--(-9.2,0)

(-1.2,0)--(0.2,0)

(1.8,0)--(3.2,0)
(3.8,0)--(5.2,0)
(5.8,0)--(7.2,0)
(7.8,0)--(9.2,0)
;

\draw (4.5,0.35) node {$v$};
\draw (1.3,-3.5) node {$X_2^3$};
\draw (1.3,3.5) node {$X_3^3$};

\draw (8.4,-1.3) node {$X^{4,4}_1$};
\draw (4.3,1.3) node {$X^{4,3}_1$};
\draw (2.3,-1.3) node {$X^{3,4}_1$};

\draw (-2.7,1.3) node {$X^{3,5}_1$};
\draw (-6.7,1.3) node {$X^{5,4}_1$};
\draw (6.7,1.3) node {$X^{4,5}_1$};

\end{tikzpicture}
\caption{\small{An $F$-partition with rewiring permutation~$\pi$. 
The permutation~$\pi$ acts on the vertices in the coloured clusters (denoted by $Y$).
Every vertex has an upper and a lower colour (the vertex $v$ has upper colour green and lower colour blue).
The lower colour of $v$ coincides with the upper colour of $\pi(v)$;
in the figure, $\pi(v)$ lies in $X_1^{3,3} \cup X_1^{3,4} \cup X_1^{3,5}$ (cf.~\eqref{rewiring:images},\eqref{rewiring:preimages}).
}}\label{fig:gadget}
\end{figure}
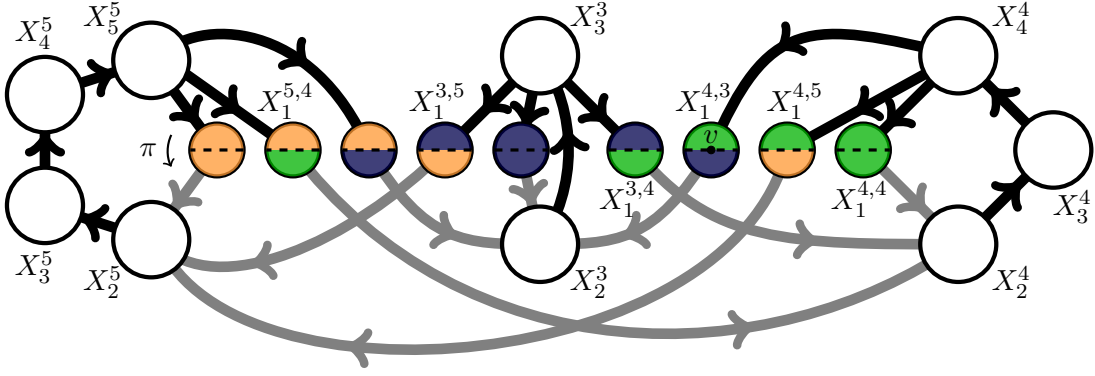

We now describe the structure of the absorbing graph, which we call an $F$-partition.
Roughly speaking, an $F$-partition $(\cX,R)$ consists of a suitable partition $\cX$ of an $n$-set~$V$ into $18$~sets and an (oriented) reduced graph~$R$ on~$\cX$. Depending on the different parts of the proof, we will endow such a partition with additional structure. For instance, the sizes of the vertex classes in $\cX$ are chosen such that there is a natural embedding of $F$ into $V$ with reduced graph~$R$ (cf.~Definition~\ref{def:nat hom}), hence the name $F$-partition. Another crucial ingredient is a `rewiring' permutation~$\pi$ on a special subset $Y\In V$ (cf.~Definition~\ref{def:rewire} and~Figure~\ref{fig:gadget}).
The absorbing graph will be a graph $G$ with vertex partition $\cX$ and reduced graph $R$ (cf.~Section~\ref{subsec:rob dec}).

\begin{defin}[$F$-partition] \label{def:gadget}
Let $F$ be a $2$-regular graph and assume that $V$ is a (vertex) set of size~$|F|$. An \defn{$F$-partition $(\cX,R)$ of $V$} is defined as follows:

\textbullet\quad
For all $a\in\Set{3,4,5}$ and $i\in [a]$, let $X^a_{i}$ be a subset of $V$ of size $\fc^F(a)$ such that 
$$\hat{\cX}:=(X^a_i)_{a\in \Set{3,4,5},i\in[a]}$$ 
is a partition of~$V$. 
This is possible by Proposition~\ref{prop:gadget sizes}\ref{prop:gadget sizes:sum}. Moreover, for each $a\in\Set{3,4,5}$, we partition $X^a_1$ further using Proposition~\ref{prop:gadget sizes}\ref{prop:gadget sizes:subpartition}, that is, let 
\begin{align}
	X^a_1=X_1^{a,3}\cupdot X_1^{a,4} \cupdot X^{a,5}_1     \label{refined partition}
\end{align}
 be a partition into sets of size $\fc^F(a,3)$, $\fc^F(a,4)$ and $\fc^F(a,5)$, respectively.
Let
\begin{align*}
	\cX & := \bigcup_{a\in\Set{3,4,5}} \Set{X_1^{a,3},X_1^{a,4},X_1^{a,5}, X^a_2,\dots,X^a_a} 
\end{align*}

\textbullet\quad  We define two (oriented) reduced graphs, one for the partition $\hat{\cX}$ and one for the refined partition $\cX$. Define $\hat{R}$ on $\hat{\cX}$ as the union of the (oriented) cycles $X^{a}_1 X^a_2 \dots X^a_a X^{a}_1$ for each $a\in\Set{3,4,5}$.
Define $R$ on $\cX$ as the union of an (oriented) path $ X^a_2 \dots X^a_a$ for each $a\in\Set{3,4,5}$ as well as the (oriented) paths $X^a_a X_1^{a,b} X^{b}_2$ for all $(a,b)\in\cI$ (cf.~Figure~\ref{fig:gadget}).

\textbullet\quad Note that $\hat{\cX}$ and $\hat{R}$ are uniquely determined by $\cX$ and $R$. Moreover, we let 
$$Y:=\bigcup_{(a,b)\in \cI}X_1^{a,b} = \bigcup_{a\in\Set{3,4,5}} X_1^a. $$
\end{defin}

Note that if at least $\eta |F|$ vertices of $F$ lie in cycles of length at least $500$, then Proposition~\ref{prop:gadget sizes}\ref{prop:gadget sizes:rich} implies that $|X|\ge \eta |F|/200$ for all $X\in \cX$.

\begin{defin}[$F$-homomorphism] \label{def:nat hom}
Let $(\cX,R)$ be an $F$-partition. An \defn{$F$-homomorphism $\sigma\colon F\to R$ for $(\cX,R)$} is defined as follows: Consider $C\in \cC(F)$. We define a homomorphism $\sigma_C\colon C\to R$ simply by walking in $R$ as indicated by~$\fa^C$. More precisely, let $\fa^C=(a_1,\dots,a_t)$ and write $$C=v_{1,1}v_{1,2}\dots v_{1,a_1}v_{2,1}\dots v_{2,a_2} \dots v_{t,1}\dots v_{t,a_t}v_{1,1}.$$ (Recall that $\sum_{i\in[t]}a_t=|C|$.) For $i\in[t]$ and $j\in[a_i]\sm\Set{1}$, define $\sigma_C(v_{i,1}):= X^{a_{i-1},a_i}_1$ (where $a_0:=a_t$) and $\sigma_C(v_{i,j}):= X^{a_i}_j$.
Let $\sigma:=\bigcup_{C\in \cC(F)}\sigma_C$. Clearly, $\sigma \colon F\to R$ is a homomorphism.
\end{defin}

We record some easy properties of such a homomorphism.

\begin{fact}\label{fact:nat hom}
Let $(\cX,R)$ be an $F$-partition and $\sigma\colon F\to R$ an $F$-homomorphism for~$(\cX,R)$. Let $\vv{F}$ denote the orientation of $F$ obtained by orienting $xy\in E(F)$ with the orientation of $\sigma(x)\sigma(y)$ in~$R$. Then $|\sigma^{-1}(X)|=|X|$ for all $X\in \cX$, and $\vv{F}$ is $1$-regular.
\end{fact}

\proof
Let $\sigma$ be defined as in Definition~\ref{def:nat hom}. Consider $C\in \cC(F)$.
Note that we have $|\sigma_C^{-1}(X^a_i)|=\fc^C(a)$ for all $a\in \Set{3,4,5}$ and $i\in[a]\sm\Set{1}$, and $|\sigma_C^{-1}(X^{a,b}_1)|=\fc^C(a,b)$ for all $(a,b)\in \cI$.
Moreover, if $\vv{C}$ denotes the orientation of $C$ obtained by orienting $xy\in E(C)$ with the orientation of $\sigma_C(x)\sigma_C(y)$ in~$R$, then $\vv{C}$ is $1$-regular.
Hence, $\vv{F}$ is clearly $1$-regular.
Moreover, for all $a\in \Set{3,4,5}$ and $i\in[a]\sm\Set{1}$, we have 
\begin{align*}
	|\sigma^{-1}(X^a_i)|=\sum_{C\in \cC(F)}|\sigma_C^{-1}(X^a_i)|=\sum_{C\in \cC(F)}\fc^C(a) \overset{\eqref{appearance counting factor}}{=} \fc^F(a) = |X^a_i|,
\end{align*}
and similarly, for all $(a,b)\in \cI$, we have $|\sigma^{-1}(X^{a,b}_1)|= \fc^F(a,b) = |X^{a,b}_1| $.
\endproof

Another property of an $F$-homomorphism $\sigma$, which will be crucial in Section~\ref{subsec:matchingII}, is that there are many subpaths of $F$ whose $\sigma$-image winds at least $5$~times around a given cycle in~$R$. This follows from Proposition~\ref{prop:gadget sizes}\ref{prop:gadget sizes:rich}.

\subsection{The robust decomposition lemma} \label{subsec:rob dec}

In this subsection, we show that given an $F$-partition~$(\cX,R)$, we can find a graph $G$ with vertex partition $\cX$ and reduced graph $R$, such that for any sparse subgraph $L\In G$, the remainder $G-L$ has an $F$-decomposition, subject to some divisibility conditions (cf.~Lemma~\ref{lem:absorbing gadget}).
This is arguably the core of our proof.

The construction of $G$ is based on a `rewiring permutation'. This permutation controls how we `merge' resolvable $C_3,C_4,C_5$-decompositions of suitable auxiliary graphs into an $F$-decomposition.

\begin{defin}[$F$-rewiring] \label{def:rewire}
Let $V$ be a set and assume that $(\cX,R)$ is an $F$-partition of~$V$. Let $\hat{\cX},\hat{R},Y$ be as in Definition~\ref{def:gadget}. An \defn{$F$-rewiring $(\pi,\pi^\ast)$ for $(\cX,R)$} is defined as follows:

\textbullet\quad We define a permutation $\pi$ on $Y$ 
that is the disjoint union of permutations $\{\pi_C\}_{C\in \cC(F)}$. 
To this end, for each $C\in \cC(F)$, let $Y_C\In Y$ be a set which consists of $\fc^C(a,b)$ vertices from $X_1^{a,b}$ for each $(a,b)\in \cI$ such that $(Y_C)_{C\in \cC(F)}$ is a partition of~$Y$. This is possible since $|X_1^{a,b}|=\fc^F(a,b)\overset{\eqref{appearance counting factor}}{=}\sum_{C\in \cC(F)}\fc^C(a,b)$ for each $(a,b)\in \cI$. 

Consider a cycle $C\in \cC(F)$ and let $\mathfrak{a}^C=(a_1,\dots,a_t)$.
For each $i\in[t]$, choose a vertex $v_i\in Y_C\cap X_1^{a_i,a_{i+1}}$ (where $a_{t+1}:=a_1$) such that the vertices $v_1,\dots,v_t$ are distinct. This is possible since $|Y_C\cap X_1^{a,b}|=\fc^C(a,b)$ for each $(a,b)\in \cI$. In particular, $Y_C=\Set{v_1,\dots,v_t}$.
Now, define $\pi_C(v_i):=v_{i+1}$ for each $i\in [t-1]$, and $\pi_C(v_t):=v_1$. (In particular, if $t=1$, then $v_1$ becomes a fixed point of~$\pi_C$.)
Finally, let $\pi\colon Y\to Y$ be the permutation which consists of all the cycles $(\pi_C)_{C\in \cC(F)}$.

Clearly, we have
\begin{align}
	\pi(X_1^{a,b}) &\In X_1^{b,3} \cup X_1^{b,4} \cup X_1^{b,5} = X_1^b,   \label{rewiring:images}   \\
	\pi^{-1}(X_1^{a,b}) &\In X_1^{3,a} \cup X_1^{4,a} \cup X_1^{5,a}     \label{rewiring:preimages}
\end{align}
for all $(a,b)\in \cI$ (cf.~Figure~\ref{fig:gadget}).

\textbullet\quad Using the permutation~$\pi$, we define a bijection $\pi^\ast$ which `rewires' edges between $Y$ and $X_2^3 \cup X_2^4\cup X_2^5$.
Let $E^{\infty}$ be the set of all edges $vv'$ such that $v\in X$, $v'\in X'$ with $XX'\in E(R)$. Similarly, let $E^{\circ\:\! \circ}$ be the set of all edges $vv'$ such that $v\in X$, $v'\in X'$ with $XX'\in E(\hat{R})$.
For an edge $e=yv\in E^{\infty}$ with $y\in Y$ and $v\in X_2^a$ for some $a\in\Set{3,4,5}$, we define
\begin{align}
	\pi^\ast(e) := \pi(y)v    \label{rewiring}
\end{align}
and $\pi^\ast(e):=e$ otherwise.
By~\eqref{rewiring:images}, we have $\pi^\ast(e)\in E^{\circ\:\! \circ}$, thus $\pi^\ast\colon E^{\infty} \to E^{\circ\:\! \circ}$.
Clearly, $\pi^\ast$ is bijective, with $\pi^{\ast -1}(e)=\pi^{-1}(y)v$ if $e=yv$ with $y\in Y$ and $v\in X_2^a$ for some $a\in\Set{3,4,5}$, and $\pi^{\ast-1}(e)=e$ otherwise.

We use $\pi^\ast$ to switch between graphs with reduced graphs $R$ and $\hat{R}$, respectively. Let $\cG^{\infty}$ be the set of all graphs $G$ on $V$ with $E(G)\In E^{\infty}$, and let $\cG^{\circ\:\! \circ}$ be the set of all graphs $G$ on $V$ with $E(G)\In E^{\circ\:\! \circ}$. Clearly, $\pi^\ast$ induces a bijection between $\cG^{\infty}$ and  $\cG^{\circ\:\! \circ}$, which we call $\pi^\ast$ again.
\end{defin}

The following are the crucial properties of our rewiring procedure.
\begin{lemma}    \label{lem:rewiring}
Let $V,\cX,R,\hat{\cX},\hat{R},Y,\pi,\pi^\ast,\cG^{\infty},\cG^{\circ\:\! \circ}$ be as in Definitions~\ref{def:gadget} and~\ref{def:rewire}.
Then the following hold:
\begin{enumerate}[label=\rm{(\roman*)}]
\item $G\in \cG^{\infty}$ is $r$-balanced with respect to $(\cX,R)$ if and only if $\pi^\ast(G)\in \cG^{\circ\:\! \circ}$ is $r$-balanced with respect to $(\hat{\cX},\hat{R})$. \label{lem:rewiring:balanced}
\item If $H\in \cG^{\circ\:\! \circ}$ consists of an $a$-partite $C_a$-factor on $(X^a_1,\dots,X^a_a)$ for each $a\in \Set{3,4,5}$, then $\pi^{\ast-1}(H)\cong F$.  \label{lem:rewiring:cycle2F}
\end{enumerate}
\end{lemma}

\proof
\ref{lem:rewiring:balanced} This follows easily from the definitions.\COMMENT{For the sake of readability, we write $G$ for the orientation $G_R$ of $G$ and $\hat{G}$ for the oriented graph $\pi^\ast(G)_{\hat{R}}$.
Clearly, we have $N^+_{G}(v)=N^+_{\hat{G}}(v)$ for all $v\in X_i^a$ with $a\in\Set{3,4,5}$ and $i\in [a]\sm \Set{1}$. Similarly, we have $N^-_{G}(v)=N^-_{\hat{G}}(v)$ for all $v\in X_i^a$ with $a\in\Set{3,4,5}$ and $i\in [a]\sm\Set{2}$.
Moreover, we have $N^+_{G}(y)=N^+_{\hat{G}}(\pi(y))$ for all $y\in Y$.
Finally, let $v \in X_2^a$ for some $a\in\Set{3,4,5}$. We have that $N^-_{G}(v),N^-_{\hat{G}}(v)\In Y$. Moreover, for $y\in Y$, we have that $y\in N^-_{G}(v)$ if and only if $\pi(y)\in N^-_{\hat{G}}(v)$. Thus, $|N^-_{G}(v)|=|N^-_{\hat{G}}(v)|$.
This completes the proof of~\ref{lem:rewiring:balanced}.}

\ref{lem:rewiring:cycle2F}
Suppose that $H\in \cG^{\circ\:\! \circ}$ consists of an  $a$-partite $C_a$-factor on $(X^a_1,\dots,X^a_a)$ for each $a\in \Set{3,4,5}$. For a vertex $y\in X^a_1$, let $C^y=yx_2^{(y)}\!\!\dots x_a^{(y)}y$ be the copy of $C_a$ which contains~$y$, where $x_i^{(y)}\in X_i^a$ for all $i\in[a]\sm\Set{1}$.
Recall that $\pi\colon Y \to Y$ is the composition of all cycles $(\pi_C)_{C\in \cC(F)}$, where $\pi_C$ is a permutation on $Y_C$ and $(Y_C)_{C\in \cC(F)}$ is a partition of~$Y$.
Clearly,
\begin{align}
	E(H) 
	&= \bigcup_{y\in Y} E(C^y)
	= \bigcup_{C\in \cC(F)}\bigcup_{y\in Y_C} E(C^y).   \label{rewiring:cycle split}
\end{align}
The crucial observation is that $\pi^{\ast-1}$ merges the cycles $(C^y)_{y\in Y_C}$ to one copy of~$C$. 

\begin{claim}
For each $C\in \cC(F)$, let $H^C$ be the graph with vertex set $\bigcup_{y\in Y_C} V(C^y)$ and edge set $\pi^{\ast-1} \left(\bigcup_{y\in Y_C} E(C^y)\right)$. Then $H^C  \cong C$.
\end{claim}

\claimproof
Fix some $C\in \cC(F)$.
Note first that for every $y\in Y_C$, the set of vertices spanned by the edges of $\pi^{\ast-1}(E(C^y))$ is $V(C^y) \cup \Set{\pi^{-1}(y)}$. Indeed, since $\pi(Y_C)=Y_C$, we have $\pi^{-1}(y)\in V(C^{y'})$ for some $y'\in Y_C$. Thus, $H^C$ is well-defined.

Now, let $\mathfrak{a}^C=(a_1,\dots,a_t)$ and write $Y_C=\Set{y_1,\dots,y_t}$, where $\pi(y_i)=y_{i+1}$ and $y_i\in X_1^{a_i,a_{i+1}}$ for all $i\in[t]$ (where $y_{t+1}:=y_1$ and $a_{t+1}:=a_1$). Consider $i\in[t]$.
Recall that $x_2^{(y_i)}$ is the neighbour of $y_i$ on $C^{y_i}$ in $X_2^{a_i}$.
By definition of $\pi^\ast$, we have that $$\pi^{\ast-1}(E(C^{y_i})) = (E(C^{y_i})-\Set{y_ix_2^{(y_i)}})\cup \Set{\pi^{-1}(y_i)x_2^{(y_i)}}.$$
Thus, $P_i:=(V(C^{y_i})\cup \Set{y_{i-1}},\pi^{\ast-1}(E(C^{y_i})))$ is a path from $y_i$ to $y_{i-1}$ of length~$a_i$. Since the internal vertices of the paths $P_1,\dots,P_t$ are mutually disjoint, we conclude that the union $H^C$ of these $t$ paths is a cycle of length $a_1+\dots+a_t=|C|$.
\endclaimproof

Since $Y_C\cap  Y_{C'}=\emptyset$ for distinct $C,C'\in \cC(F)$ and $V(C^y)\cap V(C^{y'})=\emptyset$ for distinct $y,y'\in Y$, the graphs $(H^C)_{C\in \cC(F)}$ are pairwise vertex-disjoint. Therefore, $$\pi^{\ast-1}(H) \overset{\eqref{rewiring:cycle split}}{=} \bigcup_{C\in \cC(F)}H^C \cong F.$$
\endproof

The following lemma will be used to define the absorbing graph.

\begin{lemma}  \label{lem: random regular bipartite}
Suppose $1/n \ll \eps \ll 1/s,\beta, d, 1/t$.
Let $V_1,\dots, V_t$ be sets of size~$n$ each (indices modulo~$t$). Let $\cU$ be a collection of at most $\eul^{n^{1/10}}$ sets such that for every $U\in \cU$, we have $|U|\ge \beta n$ and $U\In V_i$ for some $i\in[t]$. 
Then there exists a graph $G$ with vertex partition $(V_1,\dots,V_t)$ satisfying the following:
\begin{enumerate}[label=\rm{(\roman*)}]
\item for all $i\in [t]$, all $S\In V_{i-1}\cup V_{i+1}$ with $|S|\le s$ and each $U\in \cU$ with $U\In V_i$, we have $|U \cap \bigcap_{v\in S}N_G(v)|=(1\pm \eps)d^{|S|}|U|$;\label{lem: random regular bipartite:typical}
\item $G[V_i,V_{i+1}]$ is $\lfloor dn \rfloor$-regular for all $i\in[t]$. \label{lem: random regular bipartite:regular}
\end{enumerate}
\end{lemma}

\proof
For each $i\in[t]$, let $\cM_i$ be a decomposition of the complete bipartite graph between $V_i$ and $V_{i+1}$ into perfect matchings. Equivalently, consider an $n$-edge-colouring. Let $G$ be the (random) graph obtained by activating each $M\in \bigcup_{i\in[t]}\cM_i$ independently with probability~$d$.
For each $i\in[t]$, let $\cM_i'\In \cM_i$ be the (random) set of activated matchings. Clearly, $G[V_i,V_{i+1}]$ is $|\cM_i'|$-regular.

Using Lemma~\ref{lem:separable chernoff} and a union bound, it is easy to see that $G$ satisfies \ref{lem: random regular bipartite:typical} with probability at least $1-\eul^{-n^{1/7}}$. Moreover, standard properties of the binomial distribution yield that $G$ satisfies \ref{lem: random regular bipartite:regular} with probability at least $\Omega_d(n^{-t/2})$. Thus, there exists a graph $G$ with the desired properties.
\endproof

Later on, we will also need the following related result. 
\begin{prop}  \label{prop: random regular}
Suppose $1/n \ll \eps \ll d$. Let $V$ be a vertex set of size $n$. There exists a $2\lfloor dn/2 \rfloor$-regular $(\eps,d)$-quasirandom graph on~$V$.
\end{prop}

\proof
Let $r:=2\lfloor dn/2 \rfloor$. If $n$ is even, $K_n$ has a $1$-factorization and we can generate the desired graph as the union of random perfect matchings as above. If $n$ is odd, we take a quasirandom $r$-regular graph $G'$ on $n-1$ vertices first. Then choose a random subset $N\In V(G')$ of size~$r$. It is well known that $G'[N]$ has with high probability a perfect matching $M$ (since $r$ is even). Then removing $M$ and adding a new vertex with neighbourhood $N$ yields the desired graph.
\endproof

We are now ready to prove the robust decomposition lemma. It guarantees the existence of an absorbing graph~$G$ which not only has an $F$-decomposition itself, but even after removing a sparse balanced subgraph, it is still $F$-decomposable.

\begin{lemma} \label{lem:absorbing gadget}
Suppose $1/n \ll \eps \ll \eta,\alpha$. Let $F$ be a $2$-regular graph on $n$ vertices and assume that at least $\eta n$ vertices lie in cycles of length at least~$500$. Let $V$ be a vertex set of size~$n$ and let $(\cX,R)$ be an $F$-partition of~$V$. There exists a graph $G$ with vertex partition $\cX$ and reduced graph $R$ which satisfies the following properties:
\begin{enumerate}[label=\rm{(\roman*)}]
\item $G$ is $2r$-regular for some $r\le \alpha n$;  \label{lem:absorbing gadget:regular}
\item for all $XX'\in E(R)$, the pair $G[X,X']$ is $(\eps,d_{XX'})$-quasirandom for some $\alpha\eta/200 \le d_{XX'} \le \alpha $; \label{lem:absorbing gadget:quasirandom}
\item for every subgraph $L\In G$ which is $r_L$-balanced with respect to $(\cX,R)$ for some $r_L\le \sqrt{\eps} n$, the remainder $G-L$ has an $F$-decomposition.   \label{lem:absorbing gadget:absorb}
\end{enumerate}
\end{lemma}

\proof
We choose a new constant $s\in \N$ such that $1/n \ll \eps \ll 1/s \ll \eta,\alpha$.
Let $c:=\min\Set{\fc^F(3),\fc^F(4),\fc^F(5)}$. By Proposition~\ref{prop:gadget sizes}\ref{prop:gadget sizes:rich}, we have $c\ge \eta n/200$. Define $d_a:=\frac{\alpha c}{\fc^F(a)}$ for each $a\in\Set{3,4,5}$. Thus, $\alpha \eta/200 \le d_3,d_4,d_5 \le \alpha$. Let $r:=\lfloor \alpha c \rfloor$. Clearly, we have $\lfloor d_a \fc^F(a)\rfloor =r$ for each $a\in \Set{3,4,5}$.
Let $\hat{\cX},\hat{R}$ be as in Definition~\ref{def:gadget} and let $(\pi,\pi^\ast)$ be an $F$-rewiring for $(\cX,R)$.
For $a\in\Set{3,4,5}$, define
\begin{align*}
\cU_a &:=\Set{X_1^a,\dots,X^a_a, X_1^{a,3},X_1^{a,4},X_1^{a,5},\pi(X_1^{3,a}),\pi(X_1^{4,a}),\pi(X_1^{5,a})}.
\end{align*}
Recall from~\eqref{refined partition},~\eqref{rewiring:images} and Proposition~\ref{prop:gadget sizes}\ref{prop:gadget sizes:rich} that $X_1^{a,3},X_1^{a,4},X_1^{a,5},\pi(X_1^{3,a}),\pi(X_1^{4,a}),\pi(X_1^{5,a})$ are all subsets of $X_1^a$ of size at least $\eta n/200\ge \frac{\eta}{200} \fc^F(a)$. 

We now apply Lemma~\ref{lem: random regular bipartite} for each $a\in\Set{3,4,5}$ (with $\fc^F(a)$, $d_a$, $\eta/200$, $\cU_a$ playing the roles of $n,d,\beta,\cU$) to obtain a graph $G_a$ with vertex partition $(X_1^a,\dots,X_a^a)$ whose reduced graph is the cycle $X_1^a X_2^a \dots X^a_a X_1^a$ and which satisfies the following (indices modulo $a$):
\begin{enumerate}[label=\rm{(\alph*)}]
\item for all $i\in [a]$, all $S\In X^a_{i-1}\cup X^a_{i+1}$ with $|S|\le s$ and each $U\in \cU_a$ with $U\In X^a_i$, we have $|U \cap \bigcap_{v\in S}N_{G_a}(v)|=(1\pm \eps)d_a^{|S|}|U|$; \label{lem: random regular bipartite:typical:applied}
\item $G_a[X^a_i,X^a_{i+1}]$ is $r$-regular for all $i\in[a]$. \label{lem: random regular bipartite:regular:applied}
\end{enumerate}
For each $a\in\Set{3,4,5}$, define the density $(a\times a)$-matrix $D^a$ as $D^a_{i,i+1}:=D^a_{i+1,i}:=d_a$ for all $i\in[a]$ and $D^a_{i,i'}:=0$ otherwise. It follows immediately from \ref{lem: random regular bipartite:typical:applied} that
\begin{align}
	\mbox{$G_a$ is $(\eps,s,D^a)$-typical}.\label{blowup typical}
\end{align}
Now, define $$G:=\pi^{\ast-1}(G_3 \cup G_4 \cup G_5).$$
We claim that $G$ is the desired graph.

First, observe that $G_3 \cup G_4 \cup G_5$ is $r$-balanced with respect to $(\hat{\cX},\hat{R})$, and thus $G$ is $r$-balanced with respect to $(\cX,R)$ by Lemma~\ref{lem:rewiring}\ref{lem:rewiring:balanced}. In particular, $G$ is $2r$-regular, proving~\ref{lem:absorbing gadget:regular}.

Next, from~\ref{lem: random regular bipartite:typical:applied} and our choice of $\cU_a$, we can also deduce that for each $a\in \Set{3,4,5}$, the following pairs are $(\eps,d_a)$-quasirandom:
\begin{align*}
	\mbox{$G_a[X^a_i,X^a_{i+1}]$ for each $i\in \Set{2,\dots,a-1}$; }\\
	\mbox{$G_a[X^a_a,X_1^{a,b}]$ for each $b\in \Set{3,4,5}$; }\\
	\mbox{$G_a[\pi(X_1^{b,a}),X^a_2]$ for each $b\in \Set{3,4,5}$. }
\end{align*}
Note that for all $(a,b)\in \cI$, we have $G_a[\pi(X_1^{b,a}),X^a_2] \cong G[X_1^{b,a},X^a_2]$.\COMMENT{The isomorphism from $\pi(X_1^{b,a})\cup X^a_2$ to $X_1^{b,a}\cup X^a_2$ is the composition of $\pi^{-1}$ on $\pi(X_1^{b,a})$ and the identity on $X^a_2$}
This proves \ref{lem:absorbing gadget:quasirandom}.

It remains to prove the essential claim, which is~\ref{lem:absorbing gadget:absorb}.
To this end, let $L\In G$ be any subgraph which is $r_L$-balanced with respect to $(\cX,R)$ for some $r_L\le \sqrt{\eps} n$. Let $r':=r-r_L$.
Clearly, $G-L$ is then $r'$-balanced with respect to $(\cX,R)$. Lemma~\ref{lem:rewiring}\ref{lem:rewiring:balanced} implies that $G':=\pi^\ast(G-L)$ is $r'$-balanced with respect to $(\hat{\cX},\hat{R})$. 
For $a\in \Set{3,4,5}$, let $$G_a':=G'[X^a_1\cup \dots \cup X^a_a].$$ Hence,
\begin{align}
 G'=G_3'\cup G_4' \cup G_5'. \label{rewired split}
\end{align}
The balancedness of $G'$ implies that $d_{G_a'}(v,X^a_{i-1})=d_{G_a'}(v,X^a_{i+1})=r'$ for all $v\in X^a_i$ and $i\in[a]$ (indices modulo $a$).

Moreover, since $G_a'=G_a-\pi^\ast(L)$ and $\Delta(\pi^{\ast}(L)) \le 2\sqrt{\eps} n \le \eps^{1/3}|V(G_a)|$, Proposition~\ref{prop:typical noise} and \eqref{blowup typical} imply that $G_a'$ is still $(\eps^{1/4},s,D^a)$-typical.
Therefore, by Corollary~\ref{thm:res partite cycle dec}, $G_a'$ has a resolvable partite $C_a$-decomposition. Let $F^a_1,\dots,F^a_{r'}$ be the $C_a$-factors in such a decomposition.

Now, for each $t\in [r']$, define
\begin{align}
	F_t:=F^3_t \cup F^4_t \cup F^5_t.
\end{align}
By~\eqref{rewired split}, we have that $\set{F_t}{t\in [r']}$ is a decomposition of $G'$. Thus, $\cF:=\set{\pi^{\ast-1}(F_t)}{t\in [r']}$ is a decomposition of $\pi^{\ast-1}(G')=G-L$.
 Crucially, by Lemma~\ref{lem:rewiring}\ref{lem:rewiring:cycle2F}, $\pi^{\ast-1}(F_t)\cong F$ for all $t\in[r']$. We conclude that $\cF$ is an $F$-decomposition of~$G-L$, completing the proof.
\endproof

\subsection{Crossing edge absorption}

In this subsection, we prove the main tool for the first step of our absorption process (Lemma~\ref{lem: cross edge absorbing global}).
Roughly speaking, it states that given a graph $G$ with a vertex partition $(V_1,\ldots,V_t)$
which is quasirandom inside the clusters and sufficiently sparse between the clusters,
we can utilize copies of $F$ to cover all edges between clusters.

\begin{lemma}\label{lem: cross edge absorbing global}
Suppose $1/n \ll \eps  \ll \eta,1/t $. Let $F$ be an $n$-vertex $2$-regular graph with at least $\eta n$ vertices belonging to cycles of length at least~$15t$.
Let $(V_1,\dots, V_t)$ be a vertex partition of an $n$-vertex graph~$G$ with vertex set $V$ such that $|V_i| \geq \eta n$ and $G[V_i]$ is $(\eps ,d_i)$-quasirandom for some $d_i\ge \eta$ for all $i\in [t]$. 
Let $L$ be a graph with vertex set~$V$ such that the following hold:
\begin{enumerate}[label=\rm{(\roman*)}]
\item $\Delta(L)\le \eps n$;
\item $e_L(V_i,V_j)\ge \sqrt{\eps} n$ for all distinct $i,j \in [t]$ and $e_L(V_i)=0$ for all $i\in[t]$; \label{matching abs I cross}  \COMMENT{$1000\eps t n$ should also suffice}
\item $e_L(V_i,V\setminus V_i)$ is even for all $i\in [t]$. \label{matching abs I div}
\end{enumerate}
Then there exists a subgraph $G'$ of $G$ such that $G'\cup L$ has an $F$-decomposition.
\end{lemma}

Note that condition~\ref{matching abs I div} is clearly necessary. Moreover, we need at least a few edges between clusters (e.g.~if $F$ is a Hamilton cycle) and it turns out that the condition in \ref{matching abs I cross} suffices.

In the proof of Lemma~\ref{lem: cross edge absorbing global} we rely on Lemma~\ref{lem: cross edge absorbing}.
The statement of Lemma~\ref{lem: cross edge absorbing} is similar to the statement of Lemma~\ref{lem: cross edge absorbing global}.
However, in Lemma~\ref{lem: cross edge absorbing} $L$ is a matching which we cover with a single copy of $F$.
For the proof of Lemma~\ref{lem: cross edge absorbing} we use Lemma~\ref{lem:simple alloc}.

\begin{lemma} \label{lem:simple alloc}
Let $n,t\in \N$.
Let $F$ be an $n$-vertex $2$-regular graph where all vertices belong to cycles of length at least~$30$. Let $n_1,\dots,n_t\in \N$ be such that $\sum_{i\in[t]}n_i=n$ and $n_i\ge 50$ for all $i\in[t]$. Then there exists an assignment $f\colon V(F) \to [t]$ such that the following hold:
\begin{enumerate}[label=\rm{(\roman*)}]
\item $|f^{-1}(i)|=n_i$ for all $i\in[t]$;
\item $|f(x)-f(y)|\le 1$ for all $xy\in E(F)$;
\item the set $E$ of all edges $xy\in E(F)$ with $|\{f(x),f(y)\}|=2$ is an induced matching in $F$;
\item for all $i\in [t-1]$, there are exactly four edges $xy\in E(F)$ with $\Set{f(x),f(y)}=\Set{i,i+1}$. 
\end{enumerate} 
\end{lemma}

\proof
The proof consists of three steps.
In the first step, we define a very simple assignment $f_1:V(F)\to [t]$.
In the second step, we modify $f_1$ slightly to obtain a different assignment $f_2$.
Both assignments essentially ignore the edges of $F$ and only make sure that (i) holds with $f$ replaced by $f_1$ and $f_2$, respectively.
In the final step, it is then easy to obtain an assignment $f$ from $f_2$ which satisfies (i)--(iv).

We write $C^1,\ldots,C^s$ for the cycles in~$F$.
Let $x_1,\ldots,x_n$ be an ordering of the vertices of $F$ such that all vertices in $C^{i}$ precede all vertices in $C^{j}$ for all $i<j$.
For $i\in [t]$, let $I_i:= \{1+\sum_{j<i}n_j,2+\sum_{j<i}n_j,\ldots, \sum_{j\leq i}n_j \}$.
Clearly, $I_1,\ldots,I_t$ is a partition of $[n]$.
This also gives rise to our first assignment $f_1$ by defining $f_1(x_k):=i$ whenever $k\in I_i$.

Recall that every cycle of $F$ has length at least~$30$.
Hence for all $i\in [t]$,
there are (disjoint) intervals $I_i^-,I_i^+\sub I_i$ of size $6$
such that
\begin{enumerate}[label=(\alph*)]
	\item for all $k,k'\in I_i^-$, the vertices $x_k,x_{k'}$ belong to the same cycle of $F$ and analogously for $I_i^+$;
	\item $\max I_i^- \leq \min I_i +10$ and $\min I_i^+ \geq \max I_i -10$.
\end{enumerate}

Next, we slightly modify $f_1$. Roughly speaking, $f_2$ arises from $f_1$ by interchanging $I_i^+$ and $I_{i+1}^-$ for all $i\in [t-1]$.
To be precise, for all $k\in [n]$, let $f_2(x_k):=i$ whenever $k\in (I_i \sm (I_{i}^- \cup I_{i}^+))\cup I_{i-1}^+ \cup I_{i+1}^-$ where $I_0^+:=I_1^-$ and $I_{t+1}^-:=I_t^+$.

For all $i\in [t]$ and $j\in [s]$,
let $n(i,j):=|f_2^{-1}(i)\cap V(C^j)|$.
The following statements are easy observations that follow directly from our construction:%
\COMMENT{For (A).
The non-zero $n(i,j)$ appear consecutively in $t$. Swapping $I_i^-$ possibly creates an extra consecutive interval; $n(i-1,j)>0$ but in this case keeps $n(i,j)>0$.
If $n(i,j)>0$, note that $|(I_i \sm (I_{i}^- \cup I_{i}^+))\cup I_{i-1}^+ \cup I_{i+1}^-|\geq 20$.

For (B). This is obvious and follows from the fact that every cycle has length at least $30$.
}
\begin{enumerate}[label=(\Alph*)]
	\item for all $j\in [s]$, in the sequence $n(1,j),n(2,j),\ldots, n(t,j)$ all non-zero elements appear consecutively and if $n(i,j)>0$, then $n(i,j)\geq 6$;
	\item for all $i\in [t-1]$, there are either one or two $j\in [s]$ such that both $n(i,j)>0$ and $n(i+1,j)>0$.
\end{enumerate}
It is easy to find an assignment $f\colon V(C^j)\to [t]$ of the vertices of a single cycle $C^j$ such that
\begin{itemize}
	\item $|f^{-1}(i)\cap V(C^j)|=|f^{-1}_2(i)\cap V(C^j)|$ for all $i\in [t]$;
	\item $|f(x)-f(y)|\le 1$ for all $xy\in E(C^j)$;
	\item the set $E$ of all edges $xy\in E(C^j)$ with $|\{f(x),f(y)\}|=2$ is an induced matching in $F$;
	\item and the number of edges $xy\in E(C^j)$ with $\{f(x),f(y)\}=\{i,i+1\}$ equals $4$ if $j$ is the unique integer such that both $n(i,j),n(i+1,j)$ are positive,
$2$ if $n(i,j),n(i+1,j),n(i,j^*),n(i+1,j^*)$ are positive for some $j^*\in [s]\sm\{j\}$ and otherwise~$0$.\COMMENT{If $j$ is unique, then four crossing edges require at least six vertices (and more are ok).}
\end{itemize}
Combining these assignments for all cycles yields an assignment $f$ satisfying (i)--(iv).
\endproof

\begin{lemma}\label{lem: cross edge absorbing}
Suppose $1/n \ll \eps \ll \eta,1/t$. 
Let $F$ be an $n$-vertex $2$-regular graph with at least $\eta n$ vertices belonging to cycles of length at least~$15t$.\COMMENT{something like $\max\{3t+1,30\}$ should suffice}
Let $(V_1,\dots, V_t)$ be a vertex partition of an $n$-vertex graph~$G$ with vertex set $V$ such that $|V_i| \geq \eta n$  and
$G[V_i]$ is $(\eps ,d_{i})$-quasirandom for some $d_{i}\geq \eta$ for all $i\in [t]$. 
Let $M$ be a matching of size at most $\eps n$ such that $e_M(V_i)=0$ for all $i \in [t]$ and $e_M(V_i,V_{i+1})\ge  4$ for all $i \in [t-1]$.
Moreover, suppose that $e_M(V_i,V\setminus V_i)$ is even for all $i\in [t]$. 
Then there exists an embedding $\phi\colon F\to G\cup M$ such that $M\In \phi(F)$.
\end{lemma}

\proof
For $i\in[t]$, define $n_i:=|V_i|$. We first allocate all short cycles of $F$ to the clusters~$V_i$. To this end, let $F^<$ be the subgraph of $F$ which consists of all cycles of length less than~$15t$, and let $F^\ge$ be the subgraph of $F$ which consists of all cycles of length at least~$15t$. Let $\sigma'\colon V(F^<)\to [t]$ be such that $n_i':=n_i-|\sigma'^{-1}(i)| \ge \eta n /2t$ for all $i\in[t]$ and $\sigma'(x)=\sigma'(y)$ for all $xy\in E(F^<)$. Clearly, $\sigma'$ exists.

Our next goal is to find a function $\sigma'' \colon V(F^\ge) \to [t]$ such that the following hold:
\begin{enumerate}[label=\rm{(\alph*)}]
\item $|\sigma''^{-1}(i)|=n_i'$ for all $i\in[t]$;  \label{long cycle embed size}
\item the set $E$ of all edges $xy\in E(F^\ge)$ with $|\{\sigma''(x),\sigma''(y)\}|=2$ is an induced matching in $F$; \label{long cycle embed distance}
\item for all distinct $i,j\in [t]$, there are exactly $e_M(V_i,V_j)$ edges $xy\in E(F^\ge)$ with $\Set{\sigma''(x),\sigma''(y)}=\Set{i,j}$. \label{long cycle embed counts}
\end{enumerate}
We will find $\sigma''$ with the help of Lemma~\ref{lem:simple alloc}. Before this, we show how to complete the proof based on~\ref{long cycle embed size}--\ref{long cycle embed counts}. Assume that $\sigma''$ satisfies \ref{long cycle embed size}--\ref{long cycle embed counts}. Obviously, we will embed $E$ onto~$M$.
Let $\phi_0 \colon V(E) \to V(M)$ be such that $\phi_0(E)=M$, and $\phi_0(x) \in V_{\sigma''(x)}$ for all $x\in V(E)$. Such an embedding exists by~\ref{long cycle embed counts}. For $i\in[t]$, let $F_i:=F[\sigma'^{-1}(i)\cup \sigma''^{-1}(i)]$. Clearly, $F=\bigcup_{i\in[t]}F_i \cup E$. Now, consider $i\in[t]$. By~\ref{long cycle embed size}, we have $|F_i|=n_i$. Moreover, by~\ref{long cycle embed distance}, $V(E)\cap \sigma''^{-1}(i)$ is independent in~$F_i$. Thus, by Corollary~\ref{cor:quasirandom blowup}, there exists an embedding $\phi_i\colon F_i\to G[V_i]$ such that $\phi_i(x)=\phi_0(x)$ for all $x\in V(E)\cap \sigma''^{-1}(i)$. It is easy to see that $\phi:=\bigcup_{i\in[t]}\phi_i$ is the desired embedding, where $\phi(E)=M$.

It remains to find~$\sigma''$. Clearly, we may assume that $t\ge 2$.\COMMENT{so $15t\ge 30$} First, let $R$ be the graph on $[t]$ where $ij\in E(R)$ if and only if $e_M(V_i,V_j)$ is odd. Clearly, $R$ is Eulerian. Hence, $R$ has a decomposition $\cC$ into cycles. For $C\in \cC$ and $i\in[t]$, we define
$$g(C,i) := \begin{cases}  0  &  \mbox{if } i\notin V(C); \\
                     |C|-1   &  \mbox{if } i=\min V(C); \\
										-1  & \mbox{otherwise.}  \end{cases}$$
Note that $\sum_{i\in [t]}g(C,i)=0$ for all $C\in \cC$.
For all $i,j\in[t]$ with $|j-i|>1$, we define $\defect(i,j):=\lfloor e_M(V_i,V_j)/2 \rfloor$. Moreover, for all $i\in[t-1]$, we define $\defect(i,i+1):= \lfloor e_M(V_i,V_{i+1})/2 \rfloor-2$. By assumption, we have $\defect(i,j)\ge 0$ for all distinct $i,j\in[t]$.

Now, for all $i\in[t]$, we define
\begin{align}
	n_i^\ast := n_i' - \sum_{j=1}^{i-1} 3 \defect(i,j) + \sum_{j=i+1}^t 3 \defect(i,j)    + \sum_{C\in \cC} 3g(C,i) .  \label{new cluster sizes}
\end{align}
Observe that $\sum_{i\in[t]} n_i^\ast = \sum_{i\in[t]}n_i'=|V(F^\ge)|$. Moreover, for each $i\in[t]$, we have 
\begin{align}\label{eq:difference}
	|n_i^\ast- n_i'| 
	\leq \sum_{j\in [t]\sm \{i\}} 3 \defect(i,j) + \sum_{C\in \cC} 3|g(C,i)|
	\le 3 \eps n/2 +3t^2 \le 2\eps n.
\end{align}
\COMMENT{only need lower bound for $n_i^\ast$, but maybe helps intuition} Thus, we can apply Lemma~\ref{lem:simple alloc} to find a function $f\colon V(F^\ge) \to [t]$ such that the following hold:
\begin{enumerate}[label=\rm{(\roman*)}]
\item $|f^{-1}(i)|=n_i^\ast $ for all $i\in[t]$;  \label{simple assignment sizes}
\item $|f(x)-f(y)|\le 1$ for all $xy\in E(F^\ge)$; \label{simple assignment consecutive}
\item the set $E^\ast$ of all edges $xy\in E(F^\ge)$ with $|\{f(x),f(y)\}|=2$ is an induced matching in~$F^\ge$; \label{simple assignment ind}
\item for all $i\in [t-1]$, there are exactly four edges $xy\in E(F^\ge)$ with $\Set{f(x),f(y)}=\Set{i,i+1}$. \label{simple assignment four}
\end{enumerate}
We will obtain $\sigma''$ from $f$ by changing the image of a few vertices. Roughly speaking, we repeatedly take a subpath of $F$ which is currently embedded into $V_i$ and then move three consecutive vertices into another cluster~$V_j$. Apart from moving three vertices from $V_i$ to~$V_j$, this also produces two crossing edges (which do not share an endpoint).

More precisely, for $i\in[t]$, an \defn{$i$-target} is a subpath $P\In F^{\ge}\sm f^{-1}(V(E^\ast))$ of length~$3t$ such that $f(V(P))=\Set{i}$. Clearly, for each $i\in[t]$, there are at least $2\eps n$ vertex-disjoint $i$-targets.\COMMENT{Here we use that $n_i^\ast \ge \eta n/4t$ together with the fact that $f$ produces only few crossing edges by \ref{simple assignment consecutive} and \ref{simple assignment four}.}

For all $i\in[t]$, let $\cP_i$ be a set of $\sum_{j=i+1}^t \defect(i,j)$ vertex-disjoint $i$-targets, and let $\set{\cP_{i,j}}{j=i+1,\dots,t}$ be a partition of $\cP_i$ such that $|\cP_{i,j}|=\defect(i,j)$. For every $C\in \cC$ with $i=\min V(C)$, let $P_C$ be an $i$-target which is vertex-disjoint from all the previously chosen targets (as there are $2\eps n$ vertex-disjoint $i$-targets, $P_C$ exists by~\eqref{eq:difference}).
Let $\cP:=\bigcup_{i\in[t]}\cP_i \cup \bigcup_{C\in \cC}P_C$ be the set of all these targets.

We now define~$\sigma''$. For each target~$P\in \cP$, write $P=x_1^{(P)}\!\!\dots x_{3t+1}^{(P)}$. 
For every vertex $x\in V(F^\ge)$ which is not contained in any $P\in \cP$, we let $\sigma''(x):=f(x)$. For all $i,j\in[t]$ with $i<j$ and all $P\in \cP_{i,j}$, we define, for all $k\in[3t+1]$,
$$\sigma''(x^{(P)}_k):= \begin{cases}  j  &  \mbox{if } k\in\Set{3,4,5}; \\
										i  & \mbox{otherwise.}  \end{cases}$$
Now consider~$C\in \cC$. Write $C=i_1 i_2 \dots i_\ell i_1$ such that $i_1=\min V(C)$. We define
$\sigma''(x^{(P_C)}_{3(k-1)+s}):=i_k$ for all $k\in [\ell]\sm \Set{1}$ and $s\in \Set{0,1,2}$, and $\sigma''(x^{(P_C)}_{s}):=i_1$ for $s\in \Set{1,2}\cup\Set{3\ell,\dots,3t+1}$.

We claim that $\sigma''$ satisfies \ref{long cycle embed size}--\ref{long cycle embed counts}.
Clearly, for all $i\in[t]$, we have
\begin{align*}
	|\sigma''^{-1}(i)| & = |f^{-1}(i)| + \sum_{j=1}^{i-1} 3 |\cP_{j,i}| - \sum_{j=i+1}^t 3 |\cP_{i,j}|  - \sum_{C\in \cC} 3g(C,i)  \overset{\eqref{new cluster sizes},\ref{simple assignment sizes}}{=} n_i',
\end{align*}
thus $\sigma''$ satisfies \ref{long cycle embed size}.
Condition \ref{long cycle embed distance} clearly holds for the restriction $\sigma''{\restriction_P}$ for every $P\in \cP$.\COMMENT{As we always move $3$ vertices to the new cluster, so the two ends of the matching edges going in and out are independent in $F$} Thus, \ref{long cycle embed distance} follows from \ref{simple assignment ind} and the fact that all the paths in $\cP$ are vertex-disjoint and for every $i$-target~$P$, we have $\sigma''(x^{(P)}_s)=i$ for $s\in \Set{1,2,3t,3t+1}$.

Finally, let $i,j\in[t]$ with $i<j$. The union of all $P\in \cP_{i,j}$ gives rise to exactly $2\defect(i,j)$ edges $xy\in E(F^\ge)$ with $\Set{\sigma''(x),\sigma''(y)}=\Set{i,j}$. Moreover, if $e_M(V_i,V_j)$ is even, then there is no cycle $C\in \cC$ with $ij\in E(C)$. If $e_M(V_i,V_j)$ is odd, then there is exactly one cycle $C\in \cC$ with $ij\in E(C)$, and $P_C$ gives rise to exactly one edge $xy\in E(F^\ge)$ with $\Set{\sigma''(x),\sigma''(y)}=\Set{i,j}$.
Together with \ref{simple assignment consecutive}, \ref{simple assignment four} and the definition of $\defect(i,j)$, this implies~\ref{long cycle embed counts}.
\endproof

We can now prove Lemma~\ref{lem: cross edge absorbing global}. Essentially, we need to decompose $L$ into suitable matchings and then apply Lemma~\ref{lem: cross edge absorbing} to each matching. 

\lateproof{Lemma~\ref{lem: cross edge absorbing global}}
We call a matching $M\In L$ \defn{good} if $e(M)\le 2\eps^{1/4} n$, $e_M(V_i,V_j)\ge 4$ for all distinct $i,j\in[t]$, and $e_M(V_i,V\sm V_i)$ is even for all $i\in[t]$.
The main part of the proof is to partition $E(L)$ into good matchings. We achieve this in three steps. 
In the first step, we find a good matching $M_0$ such that for $L':=L-M_0$ the number of edges between two clusters is always even.
In the second step, we partition $L'[V_i,V_j]$ into matchings of size~$2$,
and in the final step we combine these matchings to obtain a decomposition of $L'$ into good matchings.

To this end,
let $p_{ij}:=1$ if $e_L(V_i,V_j)$ is odd and $p_{ij}:=0$ otherwise, for all distinct $i,j\in [t]$.
Let $M_0\sub L$ be a matching such that $M_0[V_i,V_j]$ consists of $4 + p_{ij}$ edges for all distinct $i,j\in [t]$.
Clearly, $M_0$ exists. Note that for each $i\in[t]$, we have $$e_{M_0}(V_i,V\sm V_i) \equiv \sum_{j\in[t]\sm\Set{i}}p_{ij} \equiv \sum_{j\in[t]\sm\Set{i}} e_L(V_i,V_j) = e_L(V_i,V\sm V_j) \equiv 0 \mod{2}.$$ Thus, $M_0$ is a good matching. Let $L':=L-M_0$. Obviously, $e_{L'}(V_i,V_j)$ is even for all distinct $i,j\in [t]$.

Next, for all distinct $i,j\in [t]$, we partition the edges of $L'[V_i,V_j]$ into matchings of size~$2$. 
Since the maximum degree of the line graph of $L'(V_i,V_j)$ is at most $2\Delta(L') \le  e_{L'}(V_i,V_j)/2 $, the complement of the line graph of $L'(V_i,V_j)$ has a perfect matching, which yields the desired partition.
We call such a matching of size~$2$ an \defn{edge pair}.

An \defn{admissible} colouring is a proper edge-colouring of $L'$ where two edges that form an edge pair receive the same colour.

We now admissibly colour $L'$ with $s:=\lceil \eps^{3/4} n \rceil$~colours. First, since $e_{L'}(V_i,V_j)- 4 s - 8t \Delta(L) \ge 2 $, it is easy to find a partial admissible colouring such that for every pair $ij\in\binom{[t]}{2}$ and each colour $k\in[s]$, exactly two edge pairs in $L'[V_i,V_j]$ are coloured~$k$.
In a second phase, we order the remaining edge pairs arbitrarily and then colour them successively, each time picking an available colour that appears least often in the current colouring. Clearly, since at each step, there are at least $s-4\Delta(L')> \eps^{3/4}n/2$ colours available for the considered edge pair, no colour class will have size more than $2\eps^{1/4}n$ in the completed colouring.\COMMENT{suppose that some colour appears $2\eps^{1/4}n$ times. In the step when the last edge was coloured in this colour, all other $s-4\Delta(L')> \eps^{3/4}n/2$ available colours must also appear at least $2\eps^{1/4}n$ times. Thus, $e(L')> 2\eps^{1/4}n \cdot \eps^{3/4}n/2$, a contradiction to $e(L')\le \eps n^2$.}

It is easy to see that each colour class of this colouring is a good matching.
To conclude, we obtained a decomposition $\cM$ of $L$ into $s+1$ good matchings. 
Finally, we apply Lemma~\ref{lem: cross edge absorbing} $s+1$ times in turn to cover all these matchings.
More precisely, we embed one copy of $F$ into the union of $M\in \cM$ and the subgraph of $G$ induced by all edges that are not covered by earlier applications of Lemma~\ref{lem: cross edge absorbing}.
Clearly, at any stage of the procedure the uncovered edges in $G[V_i]$ induce an $(\epsilon^{1/2},d_i)$-quasirandom graph for every $i\in[t]$.\COMMENT{The maximum degree of the covered part is at most $2s\le 3\eps^{3/4} n$ }
\endproof

\subsection{Atom absorption} \label{subsec:matchingII}

This subsection is devoted to the proof of the following lemma which states that if we are given a graph $G$ as in Lemma~\ref{lem:absorbing gadget}\ref{lem:absorbing gadget:quasirandom}, and a regular `leftover' inside each partition class, then we can absorb this leftover by using a few edges of~$G$. Moreover, the subgraph $A$ which we use from $G$ will be balanced with respect to $(\cX,R)$, thus ensuring that $G-A$ still has an $F$-decomposition (cf.~Lemma~\ref{lem:absorbing gadget}\ref{lem:absorbing gadget:absorb}).

\begin{lemma} \label{lem:matching abs II}
Suppose $1/n \ll \eps \ll \eta$. Let $F$ be a $2$-regular graph on $n$ vertices and assume that at least $\eta n$ vertices lie in cycles of length at least~$500$. Let $V$ be a vertex set of size~$n$ and assume that $(\cX,R)$ is an $F$-partition. Assume that $G$ is a graph with vertex partition $\cX$ and reduced graph $R$ such that for all $XX'\in E(R)$, the pair $G[X,X']$ is $(\eps,d_{XX'})$-quasirandom for some $d_{XX'}\ge \eta$. Suppose that $r\le \eps n$ and that $L_X$ is a $2r$-regular graph on $X$ for each $X\in \cX$. Let $L:=\bigcup_{X\in \cX}L_X$. Then there exists a subgraph $A\In G$ which is $r_A$-balanced with respect to $(\cX,R)$ for some $r_A\le \sqrt{\eps}n$ such that $A\cup L$ has an $F$-decomposition.
\end{lemma}

As in the previous subsection, we will prove Lemma~\ref{lem:matching abs II} by splitting $L$ into suitable matchings and then employing the blow-up lemma to extend each such matching into a copy of~$F$. 
How the latter can be done for a single matching is proved in Lemma~\ref{lem:matching abs II aux}.
Recall that since $(\cX,R)$ is an $F$-partition, Definition~\ref{def:nat hom} yields a natural homomorphism $\sigma\colon F \to R$ for which we could straightforwardly apply the blow-up lemma to find a copy of $F$ in~$G$. However, such a copy of $F$ would not contain any edges of~$L$. In order to extend a given matching into a copy of~$F$, we will locally modify~$\sigma$. Because of divisibility issues, we perform this surgery not for single matching edges, but group them into smallest balanced edge sets, which we refer to as atoms. For each atom, we will perform a slight surgery on $\sigma$ to make sure that the copy of $F$ found with the blow-up lemma will cover this particular atom.

For $a\in \Set{3,4,5}$, an \defn{$\Set{a}$-atom} is a matching of size~$a$, consisting of one edge inside $X^{a,a}_1$ and one edge inside each vertex class~$X^a_i$, for each $i\in[a]\sm\Set{1}$.
For distinct $a,b\in \Set{3,4,5}$, an \defn{$\Set{a,b}$-atom} is a matching of size~$a+b$, consisting of one edge inside each of $X^{a,b}_1$ and $X^{b,a}_1$, one edge inside each vertex class $X^a_i$, for each $i\in[a]\sm\Set{1}$, and one edge inside each vertex class $X^b_i$, for each $i\in[b]\sm\Set{1}$.
Let $\cI'$ be the set of all subsets $S$ of $\Set{3,4,5}$ with $|S|\in \Set{1,2}$.
An \defn{atom} is an $S$-atom for some $S\in \cI'$. Clearly, if $O$ is an $S$-atom, then $e(O)=\sum_{s\in S}s$.
We say that a graph\COMMENT{on $\bigcup \cX$} is \defn{internally balanced} if it is the union of edge-disjoint atoms.

\begin{fact} \label{fact:int balanced}
A graph $H$ is internally balanced if and only if the following conditions hold: 
\begin{align}
	e_{H}(X^{a}_1) = \dots = e_{H}(X^{a}_a) \mbox{ for all $a\in \Set{3,4,5}$}; \label{int balanced 1}    \\ 
	e_{H}(X^{a,b}_1) = e_{H}(X^{b,a}_1) \mbox{ for all $(a,b)\in \cI$}. \label{int balanced 2}
\end{align}
\end{fact}

\proof
Observe first that every atom satisfies~\eqref{int balanced 1} and~\eqref{int balanced 2}. Thus, if $H$ is the edge-disjoint union of atoms, it also satisfies~\eqref{int balanced 1} and~\eqref{int balanced 2}. For the converse, repeatedly remove atoms from~$H$, until this is no longer possible. Using~\eqref{int balanced 1} and~\eqref{int balanced 2}, it is easy to see that at the end of this procedure, there can be no edge left.\COMMENT{Suppose, for a contradiction, that some edge is left, say in $X^a_i$. Since $H$ satisfies~\eqref{int balanced 1} and~\eqref{int balanced 2} and we removed edge-disjoint atoms, the leftover must also satisfy~\eqref{int balanced 1} and~\eqref{int balanced 2}. Thus, by \eqref{int balanced 1}, there is a leftover edge in each of $X^a_{i'}$ for $i'\in[a]$. Let $e$ be such an edge in $X^a_1$. We distinguish two cases. If $e\in X^{a,a}_1$, then we can easily find an $\Set{a}$-atom, a contradiction. Otherwise, $e\in X^{a,b}_1$ for $b\in\Set{3,4,5}\sm \Set{a}$. Using \eqref{int balanced 2}, there must also be a leftover edge in $X^{b,a}_1$, and with~\eqref{int balanced 1} we conclude that there must be a leftover edge in each of $X^b_{i'}$ for $i'\in[b]$. This easily yields an $\Set{a,b}$-atom, a contradiction.}
\endproof

Recall that for a subgraph $A$ of $G$, where $G$ is as in the statements of Lemmas~\ref{lem:matching abs II} and~\ref{lem:matching abs II aux},
we write $A_R$ for the oriented graph obtained from $A$ by orienting the edges according to the orientation of $R$.

\begin{lemma} \label{lem:matching abs II aux}
Suppose $1/n \ll \eps \ll \eta$. Let $F$ be a $2$-regular graph on $n$ vertices and assume that at least $\eta n$ vertices lie in cycles of length at least~$500$. Let $V$ be a vertex set of size~$n$ and assume that $(\cX,R)$ is an $F$-partition. Assume that $G$ is a graph with vertex partition $\cX$ and reduced graph $R$ such that for all $XX'\in E(R)$, the pair $G[X,X']$ is $(\eps,d_{XX'})$-quasirandom for some $d_{XX'}\ge \eta$. Suppose that $\vv{M}$ is an oriented internally balanced matching of size at most $\eps n$. Then there exists a subgraph $A\In G$ such that $A\cup M  \cong F$ and for every $v\in V$, we have $d^+_{A_R\cup \vv{M}}(v) = d^-_{A_R\cup \vv{M}}(v)=1$.
\end{lemma}

\proof
Let $\cA$ be a decomposition of $M$ into atoms.
Let $\sigma\colon F\to R$ be an $F$-homomorphism for $(\cX,R)$ (cf.~Definition~\ref{def:nat hom}) and let $\vv{F}$ denote the orientation of $F$ obtained by orienting $xy\in E(F)$ with the orientation of $\sigma(x)\sigma(y)$. By Fact~\ref{fact:nat hom}, we have
\begin{align}
	\mbox{$|\sigma^{-1}(X)|=|X|$ for all $X\in \cX$} \label{cluster sizes correct}
\end{align}
and $\vv{F}$ is $1$-regular.

Observe that $R$ has a unique oriented $\ell$-cycle $D_\ell$ for each $\ell\in \Set{3,4,5,7,8,9}$, namely, for distinct $a,b\in \Set{3,4,5}$, we have 
\begin{align*}
	D_a &=X^{a,a}_1X^a_2\dots X^a_a X^{a,a}_1; \\
	D_{a+b} &= X^{b,a}_1 X^a_2\dots X^a_a X^{a,b}_1 X^b_2\dots X^b_b X^{b,a}_1.
\end{align*}

 Note that if $O$ is an atom, then $O$ consists of exactly one edge inside each vertex class $Z\in V(D_{e(O)})$.
For $\ell\in \Set{3,4,5,7,8,9}$, we say that $P\In F$ is an \defn{$\ell$-target} if $P$ is a subpath of some cycle $C\In \cC(F)$ such that $\sigma(V(P)) \In \bigcup_{X\in V(D_\ell)}X$ and $P$ has length $5 \ell$. (In other words, the $\sigma$-image of $P$ is a closed walk winding five times around the cycle~$D_\ell$.)
Our strategy is as follows.
We modify $\sigma$ (and obtain $\sigma'$) in such a way that the $\sigma'$-image of $P$ only winds four times around $D_\ell$ but `repeats' each cluster of $D_\ell$ at some point exactly once; that is, two consecutive vertices of $P$ are assigned to the same cluster.
This ensures that we can cover the edges of $O$ (which are inside the clusters) and at the same time exactly the same number of vertices are assigned to every particular cluster by $\sigma$ and $\sigma'$.

First, we assign to each atom $O\in \cA$ an $e(O)$-target $P_O$ such that all those targets are vertex-disjoint. This is possible by Proposition~\ref{prop:gadget sizes}\ref{prop:gadget sizes:rich}. Indeed, for every $a\in \Set{3,4,5}$, every appearance of $(a)^6$ in $\fa^C$ for some $C\in \cC(F)$ yields an $a$-target, and for all disinct $a,b \in \Set{3,4,5}$, every appearance of $(a,b)^6$ in $\fa^C$ for some $C\in \cC(F)$ yields an $(a+b)$-target. Hence, for every $\ell\in \Set{3,4,5,7,8,9}$, there are at least $\eta n/200$ distinct $\ell$-targets. Since $|\cA|\le |M| \le \eps n$, we can (greedily) choose an $e(O)$-target $P_O$ for each $O\in \cA$ such that all the paths $P_O$ are vertex-disjoint.

For each atom $O\in \cA$, we will now partially embed $P_O$ to cover the edges of~$O$.
Crucially, for this we use exactly as many vertices from each cluster as indicated by~$\sigma$. Thus, we will be able to apply the blow-up lemma to complete the embedding.

We now describe how the homomorphism $\sigma$ can be modified on $V(P_O)$ to absorb~$O$. Consider some $O\in \cA$ and let $\ell=e(O)$.
Write $P_O=x_1 x_2 \dots x_{5\ell+1}$. We may assume that $D_\ell=Z_1Z_2\dots Z_\ell Z_1$ is oriented as in~$R$ and $\sigma(x_1)=Z_1$. Hence, $\sigma(x_i)=Z_{i \mod{\ell}}$ for all $i\in[5\ell+1]$. Thus, $|\sigma^{-1}(Z_i) \cap V(P_O)|=5$ for $i\in [\ell]\sm\Set{1}$ and $|\sigma^{-1}(Z_1)\cap V(P_O)|=6$. Recall that $O$ consists of exactly one edge inside each vertex class $Z_i$. Hence, we can write $E(O)=\set{v_iw_i}{i\in[\ell]}$, where $v_i,w_i\in Z_i$ and $v_iw_i$ is oriented towards $w_i$ by $\vv{M}$.

We now mark $\ell$~edges of~$P_O$ which will later be embedded onto~$O$. 
More precisely, let $g:=4$ if $\ell\in\Set{4,5,7,8}$ and $g:=5$ if $\ell\in\Set{3,9}$. 
Let $$E_O:=\set{x_{2+g (i-1)}x_{3+g(i-1)}}{i\in[\ell]}.$$ Clearly, $E_O\In E(P_O)$, and
\begin{align}
	\mbox{no endvertex of $P_O$ is incident to an edge of $E_O$.} \label{fuer felix}
\end{align}
\COMMENT{$3+g(i-1)<5\ell+1$}
We now define a new function $\sigma_{O} \colon V(P_O) \to V(R)$ as follows: Let $\sigma_{O}(x_1):=Z_1 = \sigma(x_1)$. For $i\in [5\ell]$, let $\sigma_{O}(x_{i+1}):= \sigma_{O}(x_i)$ if $x_ix_{i+1}\in E_O$, and let $\sigma_{O}(x_{i+1})$ be the outneighbour of $\sigma_{O}(x_i)$ on $D_\ell$ otherwise.
For $i\in[\ell]$, let $j_i\in [\ell]$ be such that $\sigma_{O}(x_{2+g (i-1)})=Z_{j_i}$. By construction, we have $j_{i+1}\equiv j_i+(g-1) \mod{\ell}$. 
By the choice of~$g$, we have that $g-1$ and $\ell$ are coprime. Thus, $\Set{j_1,\dots,j_\ell}=[\ell]$. In other words, the $\sigma_O$-image of $P$ may be viewed as a walk in $D_\ell$ which `stops' in each $Z\in  V(D_\ell)$ exactly once.
This implies that for all $Z\in V(D_\ell)$, we have 
\begin{align}
	|\sigma_O^{-1}(Z)| = |\sigma^{-1}(Z)\cap V(P_O)|. \label{new hom sizes}
\end{align}
Note that $\sigma_{O}$ is not a homomorphism of $P_O$ into $R$. However,
\begin{align}
\mbox{$\sigma_{O}\colon P_O-E_O \to R $ is a homomorphism.} \label{new hom}
\end{align}
Moreover, it is compatible with $\sigma$ in the following sense: 
\begin{align}
	\sigma_O(x_1) = \sigma(x_1) = \sigma_O(x_{5\ell+1}) = \sigma(x_{5\ell+1}).  \label{new hom compatible}
\end{align}
Define $X_O^+:=\set{x_{2+g (i-1)}}{i\in[\ell]}$ and $X_O^-:=\set{x_{3+g (i-1)}}{i\in[\ell]}$. Moreover, define $\phi_O\colon X_O^-\cup X_O^+ \to V$ by 
\begin{align*}
	\phi_O(x_{2+g(i-1)}) &:= v_{j_i}, \\ 
	\phi_O(x_{3+g(i-1)}) &:= w_{j_i}
\end{align*}
for all $i\in[\ell]$. Hence $\phi_O(X_O^+)$ is the set of tails of $O$, and $\phi_O(X_O^-)$ is the set of heads of $O$.

Having done this for all $O\in \cA$, we define
\begin{align*}
	\vv{F'} := \vv{F}- \bigcup_{O\in \cA}E_O; \qquad X^\pm := \bigcup_{O\in \cA}X^\pm_O; \qquad  \phi' := \bigcup_{O\in \cA}\phi_O.
\end{align*}
Observe that
\begin{align}
	 d^\pm_{\vv{F'}}(x)=\begin{cases}  0 &  \mbox{if }x\in X^\pm; \\ 1   &  \mbox{otherwise.} \end{cases}   \label{in and out degrees blowup}
\end{align}
Note that since $g\ge 4$ and using~\eqref{fuer felix}, no two vertices in $X^+\cup X^-$ have a common neighbour in~$F$. Observe that $\phi'(\vv{F}[X^+\cup X^-])=\vv{M}$.
Define $\sigma' \colon V(F)\to V(R)$ by $$\sigma'(x) := \begin{cases}  \sigma_O(x) &  \mbox{if } x\in V(P_O) \mbox{ for some }O\in \cA; \\ \sigma(x)   &  \mbox{otherwise.} \end{cases}$$
By~\eqref{new hom} and~\eqref{new hom compatible}, we have that $\sigma'$ is a homomorphism from $F'$ into~$R$. Moreover, from~\eqref{cluster sizes correct} and~\eqref{new hom sizes} we can deduce that $|\sigma'^{-1}(X)| = |X|$ for all $X\in \cX$.

Finally, we can apply the blow-up lemma (Lemma~\ref{lem:blow-up new}) to obtain an embedding $\phi\colon F' \to G$ which extends $\phi'$ such that $\phi(x) \in \sigma'(x)$ for all $x\in V(F)$. Then $A:=\phi(F')$ is the desired graph. Indeed, we clearly have $A\cup M\cong F$. Moreover, from~\eqref{in and out degrees blowup} and the definition of~$\sigma'$,\COMMENT{important that $\sigma'$ still `respects' the orientation of $\sigma$} it is evident that
\begin{align*}
	 d^\pm_{A_R}(v)=\begin{cases}  0 &  \mbox{if }v\in \phi(X^\pm); \\ 1   &  \mbox{otherwise.} \end{cases}   
\end{align*}
Hence, $A_R\cup \vv{M}$ is $1$-regular, as desired.
\endproof

We will now deduce Lemma~\ref{lem:matching abs II} from Lemma~\ref{lem:matching abs II aux}. In order to partition $L$ into internally balanced matchings, it is convenient for us to exploit the following simple fact on the matching sequencibility of graphs.\COMMENT{$ms(G)$ is the maximum matching number of a linear ordering of $E(G)$, the matching number of a linear ordering is the maximal $d$ such that any $d$ consecutive edges form a matching. Introduced by Alspach}

\begin{fact} \label{fact:edge stack}
Let $G$ be a regular graph on $n$ vertices. Then there is an ordering of the edges of $G$ such that any $n/12$ consecutive edges form a matching.
\end{fact}

\proof
Suppose $G$ is $r$-regular. As the edges of $G$ can be properly coloured with $r+1$ colours, there exist matchings $M_1,\dots,M_{r+1}$ in $G$ which partition~$E(G)$ (suppose $|M_1|\le \dots \le |M_{r+1}|$). Observe that $|M_2|\ge n/4$.\COMMENT{$|M_3|+\dots+|M_{r+1}|\le (r-1)n/2 = e(G)-n/2$.} Now, assume that for some $i\in [r]$, the edges of $M_1,\dots,M_i$ are ordered such that any $n/12$ consecutive edges form a matching. We call an edge in $M_{i+1}$ \defn{blocked} if it shares a vertex with some edge in $M_i$ that belongs to the last $n/12$ edges in the current ordering. Clearly, at most $n/6$ edges of $M_{i+1}$ are blocked. Thus, we can extend the ordering by putting $n/12$ unblocked edges of $M_{i+1}$ first and then the remaining ones.\COMMENT{$|M_{i+1}|\ge n/6 + n/12$}
\endproof

We remark that the graph $A$ in Lemma~\ref{lem:matching abs II aux} is not balanced with respect to~$(\cX,R)$. However, the `inbalancedness' of $A$ is encoded in the orientation of the matching~$M$. By decomposing $L$ into suitable oriented matchings, it is not too difficult to ensure that the union of all graphs $A$ over all these matchings will indeed be balanced.

\lateproof{Lemma~\ref{lem:matching abs II}}
Note that since $L$ is $2r$-regular, there exists an orientation $\vv{L}$ of $L$ such that $\vv{L}$ is $r$-regular, that is, $d^+_{\vv{L}}(v)=d^-_{\vv{L}}(v)=r$ for all $v\in V$.

We first partition $L$ into small internally balanced matchings. For each $X\in \cX$, apply Fact~\ref{fact:edge stack} to obtain an ordering of the edges of $L_X$ such that any $\eta n/10^4$ consecutive edges form a matching (recall that $|X|\ge \eta n/200$). We view each such ordering as a stack from which we repeatedly choose the first edge which has not been chosen before. We successively choose a set of edges $M$ as follows: If for some $S\in \cI'$, there exists an $S$-atom of unchosen edges, then add those edges to~$M$ (taking the first unchosen edges from the appropriate stacks). Repeat this until $|M|\ge 2\sqrt{\eps} n$ or no atom can be found.
Clearly, $|M|\le 2\sqrt{\eps} n + 9 \le \eta n/10^4$ and thus $M$ is a matching. Moreover, since $M$ is the edge-disjoint union of atoms, it is internally balanced.
By repeating this procedure, we obtain a collection $\cM$ of edge-disjoint matchings.
Note that for all $a\in \Set{3,4,5}$ and $i\in[a]$, we have $e_L(X^a_i)=r|X^a_i|=r \fc^F(a)$. Similarly, for $(a,b)\in \cI$, using Proposition~\ref{prop:gadget sizes}\ref{prop:gadget sizes:symmetry}, we have that $e_L(X^{a,b}_1)=r|X^{a,b}_1|=r \fc^F(a,b)= r \fc^F(b,a) =r|X^{b,a}_1| = e_L(X^{b,a}_1)$. 
Clearly, since in each step we removed an atom from $L$, Fact~\ref{fact:int balanced} implies that the set of unchosen edges is the edge-disjoint union of atoms. In particular, we can find an atom as above until the very last edge of $L$ is chosen.
We conclude that $\cM$ is a decomposition of $L$, and that all matchings in $\cM$ have size at least $2\sqrt{\eps} n$, except possibly for the last one. Thus, 
\begin{align}
	|\cM| &\le \frac{e(L)}{2\sqrt{\eps} n}+1 = \frac{rn}{2\sqrt{\eps} n} +1 \le \sqrt{\eps} n.
\end{align}
From now on, we view $\cM$ as a set of oriented internally balanced matchings, where every edge simply inherits its orientation from~$\vv{L}$.

We now apply Lemma~\ref{lem:matching abs II aux} successively for each $\vv{M}\in \cM$, to find a graph $A^{M}\In G$ such that the following hold:
\begin{enumerate}[label=\rm{(\roman*)}]
\item $A^M \cup M  \cong F$;   \label{A plus match is F}
\item for every $v\in V$, we have $d^+_{A^M_R\cup \vv{M}}(v) = d^-_{A^M_R\cup \vv{M}}(v)=1$;   \label{out and in}
\item all the graphs in $\{A^M\}_{\vv{M}\in \cM}$ are pairwise edge-disjoint. \label{edge disjoint}
\end{enumerate}
 
Suppose that for some subset $\cM'\In \cM$ and all $\vv{M'}\in \cM'$ we have already found $A^{M'}$ satisfying \mbox{\ref{A plus match is F}--\ref{edge disjoint}}. 
Now we need to be able to find $A^M$. 
Let $A':=\bigcup_{\vv{M'}\in \cM'}A^{M}$ and $G':=G-A'$. 
Clearly, $\Delta(A')\le 2|\cM|\le 2\sqrt{\eps}n$, and thus $G[X,X']$ is still $(\eps^{1/3},d_{XX'})$-quasirandom whenever $XX'\in E(R)$. 
Thus, by Lemma~\ref{lem:matching abs II aux},\COMMENT{with $\eps=\eps^{1/3}$} there exists $A^M\In G'$ which satisfies \ref{A plus match is F} and \ref{out and in}. 
Clearly, $A^M$ is edge-disjoint from $A'$, thus~\ref{edge disjoint} holds as well. Hence, we can find $A^M$ for every $\vv{M}\in \cM$.

Let $A:=\bigcup_{\vv{M}\in \cM}A^M$. Clearly, $A\In G$ and $A\cup L$ has an $F$-decomposition by~\ref{A plus match is F} and~\ref{edge disjoint}. Moreover, for every $v\in V$, we have
\begin{align}
	 d^+_{A_R}(v) =  \sum_{\vv{M}\in \cM} d^+_{A^M_R}(v) \overset{\ref{out and in}}{=} |\cM| - \sum_{\vv{M}\in \cM} d^+_{\vv{M}}(v) = |\cM| - d^+_{\vv{L}}(v) = |\cM| - r.  
\end{align}
Similarly, $d^-_{A_R}(v) = |\cM|-r$ for every $v\in V$. Thus, $A$ is $r_A$-balanced with respect to $(\cX,R)$, where $r_A:=|\cM|-r$. Clearly, $r_A\le |\cM| \le \sqrt{\eps}n$.
\endproof

\subsection{Proof of Theorem~\ref{thm:main}} \label{subsec:main proof}

We are now ready to complete the proof of Theorem~\ref{thm:main}. The essential work has been done in the previous subsections. We now combine those results to finish the proof.

\lateproof{Theorem~\ref{thm:main}}
Let $F$ be as in the statement of Theorem~\ref{thm:main}.
For $\ell \in [n]$, 
let $n_\ell$ be the number of vertices in $F$ in cycles of length $\ell$. 
Of course,
\begin{align}
	\sum_{\ell=3}^n n_\ell = n.    \label{cycle length sum}
\end{align}
We may clearly assume that $\Delta\ge 2$ and that $\alpha \le \frac{1}{400\Delta}$.
Further, we choose new constants $\xi,\eps,\mu,\eta>0$ and $s\in \N$ such that
$$ 1/n \ll \xi \ll \eps \ll \mu \ll \eta \ll 1/s  \ll \alpha, 1/\Delta.$$
Let $V:=V(K_n)$.

\bigskip
\noindent {\bf Case 1: } $\sum_{\ell=3}^{500} n_\ell \ge (1-\eta)n$.
\bigskip

For the first part of the proof we assume that at least $(1-\eta) n$ vertices of $F$ are contained in cycles of length at most~$500$.
Since $F$ is obviously $\xi$-separable, by moving a suitable number of copies of $F$ from $\cF$ to~$\cH$, we may also assume that $|\cF|=\lceil \alpha n\rceil=(1\pm \xi)\alpha n$.

We will first use Theorem~\ref{thm: Padraig} to embed all graphs from~$\cH$.
Let $d:=1-2\alpha+\eps$. We partition the edges of $K_n$ into two graphs $G'$, $G''$ such that the following hold:
\begin{align}
\mbox{$G'$ is $(\xi,s,d)$-typical and $G''$ is $(\xi,s,1-d)$-typical.} \label{first splitting typical}
\end{align}
That such a partition exists can be seen easily via a probabilistic argument: for every edge independently, include it in $G'$ with probability $d$ and in $G''$ otherwise. Clearly, $K_n$ is $(\xi/2,s,1)$-typical.\COMMENT{Could also delete this and just say it follows from Proposition~\ref{prop:typical random slice}} Thus, Proposition~\ref{prop:typical random slice} implies that $G'$ and $G''$ satisfy~\eqref{first splitting typical} with high probability.

Now we use Theorem~\ref{thm: Padraig} to pack $\cH$ into~$G'$. Note that $d\ge 1-\frac{1}{200 \Delta}$. By~\eqref{first splitting typical}, we have that $d_{G'}(x)=(1\pm \xi) d n$ for all $x\in V(G')$. Moreover, $e(\cH)=\binom{n}{2}-|\cF|n \le (1-2\alpha+2\alpha\xi)\binom{n}{2}$, whereas $e(G')\ge (1-\xi)d\binom{n}{2}\ge (1-2\alpha+\eps/2)\binom{n}{2}$, implying that $e(\cH)\le (1-\eps/3)e(G')$.\COMMENT{$e(\cH)/e(G')\le \frac{1-2\alpha+2\alpha\xi}{1-2\alpha+\eps/2} = 1-\frac{\eps/2-2\alpha\xi}{1-2\alpha+\eps/2} $ and $\frac{\eps/2-2\alpha\xi}{1-2\alpha+\eps/2}\ge \frac{\eps/3}{1}$} Thus, by Theorem~\ref{thm: Padraig}, $\cH$ packs into~$G'$. Let $L_0$ be the leftover of this packing in~$G'$. Note that since $e(\cH)=\binom{n}{2}-|\cF|n \ge (1-2\alpha-2\xi)n^2/2$ and all graphs in $\cH$ are regular, at least $(1-2\alpha-2\xi)n$ edges are covered at every vertex. Thus, we have $\Delta(L_0)\le (1+\xi)dn - (1-2\alpha-2\xi)n \le 2\eps n$. We now add this leftover back to~$G''$. That is, define $G:=G''\cup L_0$.

It remains to show that $G$ has an $F$-decomposition. Using Proposition~\ref{prop:typical noise}, we can see that $G$ is $(\eps^{3/4},s,2\alpha-\eps)$-typical, say, and therefore $(\sqrt{\eps},s,2\alpha)$-typical. Note that $G$ is automatically $2|\cF|$-regular. Before starting to decompose~$G$, we first split $F$ into three subgraphs according to the cycle lengths appearing in~$F$.

Let $\ell^\ast \in \Set{3,\dots,500}$ be such that $n_{\ell^\ast} \geq n/600$. 
Clearly, $\ell^\ast$ exists since $\sum_{\ell=3}^{500} n_{\ell}\ge (1-\eta)n$.
Moreover, define
$$  I := \set{ \ell\in [500]\sm \{1,2, \ell^\ast \}}{n_\ell \geq \eta n}. $$ 
That is, $I$ consists of `significant' cycle lengths appearing in~$F$ (other than $\ell^\ast$). Note that $I$ might be empty.
Define $F_1$ to be the disjoint union of $n_{\ell^\ast}/\ell^\ast$ cycles of length $\ell^\ast$. 
Let $F_2$ be the disjoint union of $n_\ell/\ell$ cycles of length~$\ell$, for each $\ell \in I$. 
Finally, define $F_3$ as the disjoint union of $n_\ell/\ell$ cycles of length~$\ell$, for each $\ell \in [n]\sm (\Set{1,2,\ell^\ast} \cup I)$. Thus, $F$ is the disjoint union of $F_1$, $F_2$ and~$F_3$. Observe that 
\begin{align}
	|F_3|\le \eta n + 500 \eta n \le \eta^{1/2} n.   \label{bin small}
\end{align}
We will find the desired $F$-decomposition of $G$ in three steps. First, we embed all the copies of $F_2$, then those of $F_3$, and finally we complete the decomposition by embedding all copies of~$F_1$. 
In order to keep track of the vertices which we have already used to embed some part of some copy of~$F$, we introduce a set of new vertices.

Consider a vertex set $U$ of size exactly $|\cF|$ disjoint from~$V$.
Let $\tilde{G}$ be the graph on vertex set $V\cup U$ with all edges from $G$ and all edges between $V$ and~$U$. Each vertex in $U$ will represent one copy of $F$ in the required decomposition of~$G$. Let $W_F$ be the graph obtained from $F$ by adding a universal vertex. We will decompose $\tilde{G}$ into copies of $W_F$ such that the universal vertices lie in~$U$. Clearly, this corresponds to an $F$-decomposition of~$G$.

Define the $(2\times 2)$-density matrix $D$ by $D_{V}:=2\alpha$, $D_{UV}:=1$ and $D_U:=0$. Clearly, the following hold:
\begin{enumerate}[label=\rm{(\alph*)}]
\item $\tilde{G}$ is $(\sqrt{\eps},s,D)$-typical; \label{tilde G typical}
\item $d_{\tilde{G}}(v,V)=d_G(v)=2d_{\tilde{G}}(v,U)=2|\cF|$ for all $v\in V$; \label{V degree condition}
\item $d_{\tilde{G}}(u)=n$ for all $u\in U$;
\item $e_{\tilde{G}}(V)=e_{\tilde{G}}(V,U)=|\cF|n$.
\end{enumerate}
Here, \ref{tilde G typical} holds since $G$ is $(\sqrt{\eps},s,2\alpha)$-typical.\COMMENT{Similar to proof of Corollary~\ref{thm:res partite cycle dec}.}

We commence with the embedding of the copies of~$F_2$. 
To this end, partition the edge set of $\tilde{G}$ (randomly) into graphs $\Set{\tilde{G}_\ell}_{\ell\in I \cup \Set{\ell^\ast}}$, 
by including, for each $\ell \in I$, each edge independently into $\tilde{G}_\ell$ with probability $p_\ell:=\frac{n_\ell}{n}+\mu$, and into $\tilde{G}_{\ell^\ast}$ with the remaining probability $p_{\ell^\ast}:=1-\sum_{\ell\in I}p_\ell$. Note that since $\ell^\ast\notin I$, we have 
$$\sum_{\ell\in I}p_\ell=\sum_{\ell\in I}\frac{n_\ell}{n}+|I|\mu 
\overset{\eqref{cycle length sum}}{\le} 1-\frac{1}{600}+500\mu \le 1-\frac{1}{700}.$$ 
Thus, $p_{\ell^\ast}\ge \frac{1}{700}$.
Since each graph $\tilde{G}_\ell$ is a random subgraph of $\tilde{G}$, Proposition~\ref{prop:typical random slice} implies that for each $\ell\in I\cup \Set{\ell^\ast}$, with probability at least $1-1/n$, we have that
\begin{align}
	\mbox{$\tilde{G}_\ell$ is $(2\sqrt{\eps},s,p_\ell D)$-typical.} \label{tilde G l typical}
\end{align}
In particular, a decomposition of $\tilde{G}$ into graphs $\Set{\tilde{G}_\ell}_{\ell\in I \cup \Set{\ell^\ast}}$ with these properties exists. From now on, fix such a decomposition.

For $\ell\in \N\sm\Set{1,2}$, let $W_\ell$ be the wheel graph with $\ell$~spokes and hub $w$, and let $\sigma_\ell\colon V(W_{\ell})\to \Set{V,U}$ assign $w$ to $U$ and all other vertices of $V(W_\ell)$ to~$V$.

Recall that $p_\ell D_{UV}|U|=p_\ell (1\pm \xi)\alpha n= (1\pm \xi) p_\ell D_V |V|/2$. Thus, for each $\ell\in I$ separately, we can apply Corollary~\ref{cor:approximate wheel} to find a $(W_\ell,\sigma_\ell)$-packing $\cW_\ell$ in $\tilde{G}_\ell$ such that the leftover $L_\ell$ satisfies $\Delta(L_\ell) \le \mu^2 n$.\COMMENT{possible since $\eps \ll \mu$}

\begin{NoHyper}
\begin{claim}
For every $\ell \in I$, there exists $\cW_\ell'\In \cW_\ell$ such that every $u\in U$ is the hub of exactly $n_\ell/\ell$ wheels in~$\cW_\ell'$, and the new leftover $L_\ell':=\tilde{G}_\ell- \cW_\ell'$ satisfies $\Delta(L_\ell')\le 25\sqrt{\mu}n$.
\end{claim}
\end{NoHyper}

\claimproof
Fix some $\ell \in I$.
Consider $u\in U$. Let $s_u'$ be the number of wheels in $\cW_\ell$ which contain~$u$. Clearly, $d_{L_\ell}(u)=d_{\tilde{G}_\ell}(u)-s_u'\cdot \ell$. 
The typicality of $\tilde{G}_\ell$ implies that $d_{\tilde{G}_\ell}(u) = (1\pm 2\sqrt{\eps})p_\ell n =  n_\ell +\mu n \pm 2\sqrt{\eps}n $. 
Since $\Delta(L_\ell)\le \mu^2 n$, we conclude that $s_u'=(d_{\tilde{G}_\ell}(u)-d_{L_\ell}(u))/\ell = \frac{n_\ell}{\ell} +\frac{\mu n}{\ell} \pm 2\mu^2 n$. 
In particular, $s_u'\ge \frac{n_\ell}{\ell}$. Let $s_u:=s_u'-\frac{n_\ell}{\ell}$. 
Note that $s_u\le 2\mu n/\ell$.

For each $u\in U$, we want to discard $s_u$ wheels from $\cW_\ell$ which contain~$u$. We can do this independently for each $u\in U$, however we need to be careful not to delete too many wheels which contain a particular vertex $v\in V$. Thus, for each $u\in U$, we pick $s_u$ wheels from $\cW_\ell$ which contain~$u$ uniformly at random, independently of the choices for other vertices $u'\in U$.
Clearly, this yields $\cW_\ell'\In \cW_\ell$ such that every $u\in U$ is the hub of exactly $n_\ell/\ell$ wheels in~$\cW_\ell'$.
Let $L_\ell'$ be the new leftover. Clearly, $d_{L_\ell'}(u)\le \Delta(L_\ell)+\ell \cdot s_u \le 3\mu n$ for all $u\in U$. Now, consider $v\in V$.
If $W\in \cW_\ell$ contains $v$ and $u$, then the probability that $W$ is discarded is $s_u/s_u'$. 
Since $\ell\in I$, we have $s_u'\ge n_\ell/\ell \ge \eta n/\ell$ and thus $s_u/s_u'\le \sqrt{\mu}$. 
With Lemma~\ref{lem:chernoff}\ref{chernoff crude},%
\COMMENT{Since for each $u\in U$, there is only one $W\in \cW_\ell$ containing both $u$ and $v$ for all $W\in \cW_\ell$ containing $v$, the events $W\in \cW\setminus \cW'$ are independent.}
we conclude that with high probability, for all $v\in V$, at most $7 \sqrt{\mu} |U|$ wheels containing $v$ are discarded from~$\cW_\ell$, implying that $d_{L_\ell'}(v) \le \Delta(L_\ell) + 7 \sqrt{\mu} |U|\cdot 3 \le 25 \sqrt{\mu} n$.
Thus, there exists a choice of $\cW_\ell'$ for which $\Delta(L_\ell')\le 25 \sqrt{\mu} n$.
\endclaimproof

Recall that each $u\in U$ encodes a copy $F_u$ of $F$ in the desired decomposition of~$G$. At this stage, if we consider all wheels in $\bigcup_{\ell\in I}\cW_\ell'$ which contain~$u$, this yields a copy $F_{2,u}$ of~$F_2$. This concludes the first step.

\medskip
To continue, we combine all the leftover graphs $L_\ell'$ with $\tilde{G}_{\ell^\ast}$, that is, we define $$G^\ast:= \tilde{G}_{\ell^\ast}\cup \bigcup_{\ell\in I}L_\ell'.$$
By Proposition~\ref{prop:typical noise} and~\eqref{tilde G l typical}, $G^\ast$ is still $(\mu^{1/3},s,p_{\ell^\ast} D)$-typical.
In the next step, we embed all copies of~$F_3$. For each $u\in U$, we want to embed a copy $F_{3,u}$ of $F_3$ into $G^\ast[N_{G^\ast}(u)]$.
 To ensure that these embeddings are edge-disjoint, we proceed sequentially. Suppose that for some subset $U'\In U$, we have already successfully embedded such copies $F_{3,u'}$ edge-disjointly whilst ensuring that
\begin{align}
\Delta(F_{3,U'})\le 2\eta^{1/4}n,  \label{greedy embedding maxdeg}
\end{align}
where $F_{3,U'}:=\bigcup_{u'\in U'}F_{3,u'}$.
Now we want to find $F_{3,u}$ such that the above holds with $U'$ replaced by $U'\cup \Set{u}$. 
Let $B:=\set{v\in V}{d_{F_{3,U'}}(v)\ge \eta^{1/4}n}$. Since $e(F_{3,U'})\le |U| \cdot |F_3|\le \eta^{1/2}n^2$ by~\eqref{bin small}, we deduce that $|B|\le 2\eta^{1/4}n$. We exclude the vertices of $B$ and the edges of $F_{3,U'}$ when finding~$F_{3,u}$. Let $G_u:=G^\ast[N_{G^\ast}(u)\sm B]-F_{3,U'}$.
\begin{NoHyper}
\begin{claim} \label{claim:greedy embedding}
Every pair of vertices in $G_u$ has at least $p_{\ell^\ast}^3\alpha^2 n$ common neighbours (in~$G_u$).
\end{claim}
\end{NoHyper}

\claimproof
Let $v,v'$ be two vertices in $G_u$. Since $G^\ast$ is $(\mu^{1/3},s,p_{\ell^\ast} D)$-typical, we have that $v,v',u$ have $(1\pm \mu^{1/3})p_{\ell^\ast}^3 D_{UV}D_V^2 n$ common neighbours in~$G^\ast$. In particular, $v$ and $v'$ have at least $(1-\mu^{1/3}) 4\alpha^2 p_{\ell^\ast}^3 n$ common neighbours in $G^\ast[N_{G^\ast}(u)]$. At most $|B|+2\Delta (F_{3,U'}) \le 6\eta^{1/4} n$ have to be discarded.
\endclaimproof

By Claim~\ref*{claim:greedy embedding} and since $|F_3|\le \eta^{1/2}n$, we can greedily find a copy $F_{3,u}$ of $F_3$ in $G_u$ by choosing one vertex after the other. Clearly, \eqref{greedy embedding maxdeg} still holds with $U'$ replaced by $U'\cup \Set{u}$.
Thus, we can carry out this embedding for all $u\in U$, which completes the second step.

\medskip
It remains to embed all copies of $F_1$. For each $u\in U$, let $W_{3,u}$ be the graph obtained from $F_{3,u}$ by adding all edges from $V(F_{3,u})$ to~$u$. Moreover, let $\cW^3:=\set{W_{3,u}}{u\in U}$. Since $|F_3|\le \eta^{1/2}n$ and using \eqref{greedy embedding maxdeg}, we have that $\Delta(\cW^3)\le 3\eta^{1/4}n$.\COMMENT{vertices in $U$ have degree at most $|F_3|$, a vertex in $V$ is contained in at most $\eta^{1/4}n$ wheels and thus is joined to at most $\eta^{1/4}n$ vertices in $U$, which gives $3\eta^{1/4}n$.}
Finally, define $$\hat{G}:=G^\ast-\cW^3.$$
By Proposition~\ref{prop:typical noise}, $\hat{G}$ is still $(\eta^{1/5},s,p_{\ell^\ast} D)$-typical. 
Recall that $p_{\ell^\ast}\ge {1}/{700}$ and $\eta \ll 1$. 
Moreover, $\hat{G}$ is obtained from $\tilde{G}$ by removing edge-disjoint wheel graphs with hubs in~$U$. 
Thus, from~\ref{V degree condition}, we have $d_{\hat{G}}(v,V)=2d_{\hat{G}}(v,U)$ for all $v\in V$. 
We also have $d_{\hat{G}}(u)=n-|F_2|-|F_3|=|F_1|= n_{\ell^\ast}$ for all $u\in U$.
Hence, by Corollary~\ref{thm:non partite cycle dec}, $\hat{G}$ has a $(W_{\ell^\ast},\sigma_{\ell^\ast})$-decomposition~$\cW_{\ell^\ast}$. For $u\in U$, let $F_{1,u}$ be the collection of all $\ell^\ast$-cycles which together with $u$ form a wheel in~$\cW_{\ell^\ast}$. Clearly, $F_{1,u}$ is a copy of $F_1$.

Therefore, for every $u\in U$, the graph $F_{1,u}\cup F_{2,u} \cup F_{3,u}$ is a copy of $F$ in~$G$, and all these copies are edge-disjoint. Hence, $G$ has an $F$-decomposition.
This completes the proof of Case~1.

\bigskip
\noindent {\bf Case 2: } $\sum_{\ell=500}^n n_\ell \ge \eta n$.
\bigskip

We assume now that at least an $\eta$-fraction of vertices in $F$ lie in cycles of length at least~$500$.
Our strategy is as follows.
Let $(\cX,R)$ be an $F$-partition of~$V$.
We select the following three edge-disjoint graphs.
Given $(\cX,R)$, we set aside an absorbing graph $G^{abs}$ as in Lemma~\ref{lem:absorbing gadget} and 
also reserve a regular graph $G_\cX$ which is the union of quasirandom graphs $G_X$ on $X$ for each $X\in \cX$.
We also set aside a random edge slice $\tilde{G}$ of the remaining graph.
The graph $\tilde{G}$ is much sparser than $G_\cX$,
and $G_\cX$ is much sparser than $G^{abs}$.
Then we apply the bandwidth theorem for approximate decompositions (Theorem~\ref{thm: Padraig}) to find a packing of~$\cH$ in the remainder.
This yields a very sparse uncovered leftover.
Afterwards, we add $\tilde{G}$ to the leftover from this packing to make it sufficiently well behaved.
By Lemma~\ref{lem: cross edge absorbing global}, we can cover this new leftover with a few copies of~$F$ by using (additionally) only edges of~$G_\cX$.
In a further step, we utilize Lemma~\ref{lem:matching abs II}
to cover the remaining edges of~$G_\cX$ by using very few edges of~$G^{abs}$.
In particular, by Lemma~\ref{lem:absorbing gadget} the remaining subgraph of $G^{abs}$ will still admit an $F$-decomposition.

Now we turn to the details. Let $(\cX,R)$ be an $F$-partition of~$V$.
First, we apply Lemma~\ref{lem:absorbing gadget} (with $\eta$ playing the roles of $\alpha$ and $\eta$, and $\mu$ playing the role of $\eps$) to find a graph $G^{abs}$ with vertex partition $\cX$ and reduced graph $R$ which satisfies the following properties:
\begin{enumerate}[label=\rm{(\roman*)}]
\item $G^{abs}$ is $2r$-regular for some $r\le \eta n$;  \label{lem:absorbing gadget:regular:applied}
\item for all $XX'\in E(R)$, the pair $G^{abs}[X,X']$ is $(\mu,d_{XX'})$-quasirandom for some $\eta^2/200 \le d_{XX'} \le \eta$; \label{lem:absorbing gadget:quasirandom:applied}
\item for every subgraph $A\In G^{abs}$ that is $r_A$-balanced with respect to $(\cX,R)$ for some $r_A\le \sqrt{\mu} n$, the remainder $G^{abs}-A$ has an $F$-decomposition.   \label{lem:absorbing gadget:absorb:applied}
\end{enumerate}

Let $c:=\min_{(a,b)\in \cI}\fc^F(a,b)$. By Proposition~\ref{prop:gadget sizes}\ref{prop:gadget sizes:rich}, we have $c\ge \eta n/200$. 
Define $d_{X^a_i}:=\frac{\mu c}{\fc^F(a)}$ for each $a\in\Set{3,4,5}$ and $i\in [a]\sm \Set{1}$, and set $d_{X^{a,b}_1}:=\frac{\mu c}{\fc^F(a,b)}$ for all $(a,b)\in \cI$. 
Thus, $\mu \eta/200 \le d_X \le \mu$ for all $X\in \cX$. 
Let $r^{\circ}:=\lfloor \mu c/2 \rfloor$. 
Clearly, we have $\lfloor d_X |X|/2 \rfloor = r^{\circ}$ for all $X\in \cX$.
Now, for each $X\in \cX$, we apply Proposition~\ref{prop: random regular} to find a graph $G_X$ on $X$ which is $2r^{\circ}$-regular and $(\eps,d_X)$-quasirandom.
Let $$G_\cX:=\bigcup_{X\in \cX} G_X.$$

Observe that $K_n-G^{abs} -G_\cX$ is an $(n - 1 - 2r - 2r^\circ)$-regular graph.
Next, we select in $K_n-G^{abs} -G_\cX$ every edge independently with probability $\eps \cdot \frac{n-1}{n - 1 - 2r - 2r^\circ}$.
Using Lemma~\ref{lem:chernoff}\ref{chernoff eps}, with probability at least $1/2$, 
this yields a graph $\tilde{G}$
where $d_{\tilde{G}}(v)=\eps n \pm \xi n$ for all $v\in V$ and $e_{\tilde{G}}(X,X')\geq \eps^2 n^2 \ge \sqrt{2\eps}n$ for all distinct $X,X'\in \cX$.
Let $\tilde{G}$ be some graph with these properties.\COMMENT{need only linearly many edges between any two clusters. Could also find $\tilde{G}$ ad hoc. In particular, a regular graph would be nice subsequently. Hamilton cycle?}
Now, let $$G':= K_n - G^{abs} -  G_\cX - \tilde{G}. $$
We will first use Theorem~\ref{thm: Padraig} to pack $\cH$ into~$G'$. By moving some copies of $F$ from $\cF$ to $\cH$,
we may assume that $|\cF|=r + r^\circ + \lfloor\eps n\rfloor$.\COMMENT{$r + r^\circ + \lfloor\eps n\rfloor \le \alpha n$}
Let $d:=1-\frac{2r+2r^\circ}{n}-\eps$. Note that $d\ge 1-\frac{1}{200 \Delta}$ by~\ref{lem:absorbing gadget:regular:applied}.\COMMENT{With a lot of room to spare}
Clearly, $d_{G'}(v)= (d\pm 2\xi)n$ for all $v\in V$. Observe that $e(\cH)=\binom{n}{2}-|\cF|n \le (d-0.9\eps)n^2/2$ and hence $e(\cH)\le (1-\eps/2)e(G')$.\COMMENT{$e(G')\ge (d-2\xi)n^2/2$ and thus $e(\cH)/e(G')\le (d-0.9\eps)/(d-2\xi) = 1 - (0.9\eps-2\xi)/(d-2\xi)$ and $(0.9\eps-2\xi)/(d-2\xi)\ge \eps/2$.}
Thus, by Theorem~\ref{thm: Padraig}, $\cH$ packs into~$G'$. Let $L_1$ be the leftover of this packing in~$G'$, and let $L_1':=L_1 \cup \tilde{G}$. 
It remains to show that $G^{abs} \cup G_\cX  \cup L_1'$ has an $F$-decomposition.

First note that because all graphs in $\cH$ are regular, $G^{abs}$ is $2r$-regular, and $G_\cX$ is $2r^\circ$-regular, we can conclude that
\begin{align}
	\mbox{$d_{L_1'}(v) = 2|\cF| - 2r  - 2r^\circ  =  2\lfloor\eps n\rfloor$ for all $v\in V$.}  \label{leftover regular}
\end{align}
\COMMENT{$d_{L_1'}(v)= n-1 - 2r-2r^\circ - d_{\cH}$ and $n-1=d_{\cH}+ 2|\cF|$}

Let $L_1'':=L_1' - \bigcup_{X\in \cX}L_1'[X] $ be the subgraph of $L_1'$ which consists of all the `crossing' edges.

Clearly, $\Delta(L_1'') \le \Delta(L_1') \le 2\eps n$ by~\eqref{leftover regular}. Moreover, for all distinct $X,X'\in \cX$, we have $e_{L_1''}(X,X')\geq e_{\tilde{G}}(X,X') \ge \sqrt{2\eps} n$.
Crucially, for all $X\in \cX$, we have that
\begin{align*}
	e_{L_1''}(X,V\sm X) &= e_{L_1'}(X,V\sm X)  = \sum_{v\in X}d_{L_1'}(v) - 2e_{L_1'}(X)  \overset{\eqref{leftover regular}}{\equiv}  0\mod{2} .
\end{align*}
Hence, we can employ Lemma~\ref{lem: cross edge absorbing global} (with $2\eps,\mu\eta/200,18$ playing the roles of $\eps,\eta,t$) to obtain a subgraph $G_\cX'$ of $G_\cX$ such that $G_\cX'\cup L_1''$ has an $F$-decomposition.

We now define
$$L_2:= (G_\cX-G_\cX') \cup \bigcup_{X\in \cX}L_1'[X] = (G_\cX \cup L_1') - (G_\cX'\cup L_1'').$$
Note that it remains to show that $G^{abs} \cup L_2 = (G^{abs} \cup G_\cX \cup L_1') - (G_\cX'\cup L_1'') $ has an $F$-decomposition.
Clearly, both endpoints of any edge in $L_2$ lie in the same part $X\in \cX$.
Moreover, as $d_{L_2}(v) = 2r^\circ + 2\lfloor\eps n\rfloor - d_{G_\cX'\cup L_1''}(v)$ for all $v\in V$ and $G_\cX'\cup L_1''$ is even-regular\COMMENT{has an $F$-decomposition}, we deduce that $L_2$ is $2s$-regular for some $s\in \N$ with $s\leq \mu n/2 +\eps n \le \mu n$.
Hence, by~\ref{lem:absorbing gadget:quasirandom:applied}, we can apply Lemma~\ref{lem:matching abs II} (with $\mu$ playing the role of $\eps$)
to find a subgraph $A\In G^{abs}$ such that $A \cup L_2$ has an $F$-decomposition and $A$ is $r_A$-balanced with respect to~$(\cX,R)$ for some $r_A \leq\sqrt{\mu}n$.
Finally, by~\ref{lem:absorbing gadget:absorb:applied}, $G^{abs} - A$ also has an $F$-decomposition, which completes the proof.
\endproof

\section{Concluding remarks} \label{sec:the end}

Note that our main result (Theorem~\ref{thm:main}) is formulated in terms of decompositions
of the complete graph $K_n$.
Essentially the same argument allows us to replace  $K_n$ by any $n$-vertex regular host graph $G$
which is almost complete in the sense that the degrees are $n-o(n)$. This allows us to find an `Oberwolfach factorization' where the first $o(n)$ $F$-factors in $K_n$ can be chosen arbitrarily.

\begin{theorem} \label{thm:min deg version}
For given $\Delta\in \mathbb{N}$ and $\alpha >0$, there exist $\xi_0>0$ and $n_0 \in \mathbb{N}$ such that the following holds for all $n\geq n_0$ and $\xi <\xi_0$. Let $G$ be an $r$-regular $n$-vertex graph with $r\ge (1-\xi)n$ and let $\cF, \cH$ be collections of graphs satisfying the following:
\begin{itemize}
\item $\cF$ is a collection of at least $\alpha n$ copies of $F$, where $F$ is a $2$-regular $n$-vertex graph;
\item each $H\in \cH$ is a $\xi$-separable  $n$-vertex $r_H$-regular graph for some $r_H\leq \Delta$;
\item $e(\cF\cup \cH)=e(G)$.
\end{itemize}
Then $G$ decomposes into $\cF\cup \cH$.
\end{theorem}

The proof of Theorem~\ref{thm:min deg version} is essentially the same as that of Theorem~\ref{thm:main}.
The only significant difference is the construction of the absorbing graph $G$ in the proof of Lemma~\ref{lem:absorbing gadget}.
In Lemma~\ref{lem:absorbing gadget}, it was sufficient to prove simply the existence of such a graph.
We achieved this by first proving the existence of suitable `permuted blow-ups' $G_a$ for $a\in \Set{3,4,5}$.
The existence of the $G_a$ in turn followed from Lemma~\ref{lem: random regular bipartite}.
We now need to find $G$ in an $r$-regular graph $G^*$ of degree $r\ge (1-\xi)n$.
For this, we argue as in the proofs of Lemmas~\ref{lem: random regular bipartite} and~\ref{lem:absorbing gadget}:
We first choose an $F$-partition $(\cX, R)$ of $V=V(G^*)$ and let $G^{**}$ denote the spanning subgraph of  $G^*$ containing all those edges which correspond to edges of $R$.
Then for each $a\in \{3,4,5\}$, the graph $\pi^*(G^{**})[X^a_1 \cup \dots \cup X^a_a]$ is a blown-up cycle. Our aim now is to construct a graph $G_a\subseteq \pi^*(G^{**})[X^a_1 \cup \dots \cup X^a_a]$
satisfying~\ref{lem: random regular bipartite:typical:applied} and~\ref{lem: random regular bipartite:regular:applied} in the proof of Lemma~\ref{lem:absorbing gadget}. This can be done similarly as before: since for each $i\in [a]$ the pair $\pi^*(G^{**})[X^a_i,X^a_{i+1}]$ has
minimum degree at least $(1-\sqrt{\xi})|X^a_i|$, it contains at least $(1-\xi^{1/3})|X^a_i|$ edge-disjoint perfect matchings.
Choosing a suitable number of these at random as in the proof of Lemma~\ref{lem: random regular bipartite} again gives the desired graph~$G_a$. 

\COMMENT{A similar modification arises in the construction of $G_{\mathcal{X}}$ in Case~2 of the proof of Theorem~\ref{thm:main}. But regularity is not so essential here.
To construct the quasi-random graph $G_X$ on $X$ (for each $X\in\mathcal{X}$), we argue similarly as in the proof of Lemma~\ref{prop: random regular}, utilizing the fact that $G^*[X]$ is almost complete and thus can be almost decomposed into matchings of size $\lfloor |X|/2 \rfloor$.}

Finally, we remark that Keevash and Staden~\cite{KS:20a} have recently given a new proof of our result. Their approach is also based on~\cite{keevash:18b} and uses probabilistic techniques and the absorption method.
They do not require the condition $m_1\ge \alpha n$ in Theorem~\ref{thm:multicity}, that is, they solve the problem when $F_1,\dots,F_{(n-1)/2}$ can be arbitrary $2$-factors.

\bibliographystyle{amsplain_v2.0customized}
\bibliography{References}

\end{document}